\documentclass{article}%
\usepackage{amssymb}
\usepackage{amsmath}%
\usepackage{amsfonts}%
\usepackage{graphicx}
\providecommand{\U}[1]{\protect\rule{.1in}{.1in}}
\newtheorem{theorem}{Theorem}[section]

\newtheorem{corollary}[theorem]{Corollary}

\newtheorem{lemma}[theorem]{Lemma}

\newtheorem{proposition}[theorem]{Proposition}

\newenvironment{proof}[1][Proof]{\textbf{#1.} }{\ \rule{0.5em}{0.5em}}
\begin{document}

\title{Canonical Forms for Unitary Congruence and *Congruence }
\author{Roger A. Horn\thanks{Mathematics Department, University of Utah, Salt Lake
City, Utah, USA 84103, \texttt{rhorn@math.utah.edu}}\quad and Vladimir V.
Sergeichuk\thanks{Institute of Mathematics, Tereshchenkivska 3, Kiev, Ukraine,
\texttt{sergeich@imath.kiev.ua. } Partially supported by FAPESP (S\~{a}o
Paulo), processo 05/59407-6\texttt{. }}}
\date{}

\maketitle

\begin{abstract}
We use methods of the general theory of congruence and *congruence for complex
matrices---regularization and cosquares---to determine a unitary congruence
canonical form (respectively, a unitary *congruence canonical form) for
complex matrices $A$ such that $\bar{A}A$ (respectively, $A^{2}$) is normal.

As special cases of our canonical forms, we obtain---in a coherent and
systematic way---known canonical forms for conjugate normal, congruence
normal, coninvolutory, involutory, projection, $\lambda$-projection, and
unitary matrices. But we also obtain canonical forms for matrices whose
squares are Hermitian or normal, and other cases that do not seem to have been
investigated previously.

We show that the classification problems under (a) unitary *congruence when
$A^{3}$ is normal, and (b) unitary congruence when $A\bar{A}A$ is normal, are
both unitarily wild, so there is no reasonable hope that a simple solution to
them can be found.

\end{abstract}

\section{Introduction}

We use methods of the general theory of congruence and *congruence for complex
matrices---regularization and cosquares---to determine a unitary congruence
canonical form (respectively, a unitary *congruence canonical form) for
complex matrices $A$ such that $\bar{A}A$ (respectively, $A^{2}$) is normal.

We prove a regularization algorithm that reduces any singular matrix by
unitary congruence or unitary *congruence to a special block form. For
matrices of the two special types under consideration, this special block form
is a direct sum of a nonsingular matrix and a singular matrix; the singular
summand is a direct sum of a zero matrix and some canonical singular 2-by-2
blocks. Analysis of the cosquare and *cosquare of the nonsingular direct
summand reveals 1-by-1 and 2-by-2 nonsingular canonical blocks.

As special cases of our canonical forms, we obtain---in a coherent and
systematic way---known canonical forms for conjugate normal, congruence
normal, coninvolutory, involutory, projection, and unitary matrices. But we
also obtain canonical forms for matrices whose squares are Hermitian or
normal, $\lambda$-projections, and other cases that do not seem to have been
investigated previously. Moreover, the meaning of the parameters in the
various canonical forms is revealed, along with an understanding of when two
matrices in a given type are in the same equivalence class.

Finally, we show that the classification problems under (a) unitary
*congruence when $A^{3}$ is normal, and (b) unitary congruence when $A\bar
{A}A$ is normal, are both unitarily wild, so there is no reasonable hope that
a simple solution to them can be found.

\section{Notation and definitions}

All the matrices that we consider are complex. We denote the set of $n$-by-$n
$ complex matrices by $M_{n}$. The \textit{transpose} of $A=[a_{ij}]\in M_{n}
$ is $A^{T}=[a_{ji}]$ and the \textit{conjugate transpose} is $A^{\ast}%
=\bar{A}^{T}=[\bar{a}_{ji}]$; the \emph{trace} of $A$ is $\operatorname{tr}%
A=a_{11}+\cdots+a_{nn}$.

We say that $A\in M_{n}$ is: \emph{unitary} if $A^{\ast}A=I$;
\emph{coninvolutory} if $\bar{A}A=I$; a $\lambda$-\emph{projection} if
$A^{2}=\lambda A$ for some $\lambda\in\mathbb{C}$ (\emph{involutory} if
$\lambda=1$); \emph{normal} if $A^{\ast}A=AA^{\ast}$; \emph{conjugate normal}
if $A^{\ast}A=\overline{AA^{\ast}}$; \emph{squared normal} if $A^{2}$ is
normal; and \emph{congruence normal} if $\bar{A}A$ is normal. For example, a
unitary matrix is both normal and conjugate normal; a Hermitian matrix is
normal but need not be conjugate normal; a symmetric matrix is conjugate
normal but need not be normal.

If $A$ is nonsingular, it is convenient to write $A^{-T}=(A^{-1})^{T}$ and
$A^{-\ast}=(A^{-1})^{\ast}$; the \emph{cosquare} of $A$ is $A^{-T}A$ and the
\emph{*cosquare} is $A^{-\ast}A$.

We consider the \emph{congruence} equivalence relation ($A=SBS^{T}$ for some
nonsingular $S$) and the finer equivalence relation \emph{unitary congruence}
($A=UBU^{T}$ for some unitary $U$). We also consider the \emph{*congruence}
equivalence relation ($A=SBS^{\ast}$ for some nonsingular $S$) and the finer
equivalence relation \emph{unitary *congruence} ($A=UBU^{\ast}$ for some
unitary $U$). Two pairs of square matrices of the same size $(A,B)$ and
$(C,D)$ are said to be \emph{congruent}, and we write $(A,B)=S(C,D)S^{T}$, if
there is a nonsingular $S$ such that $A=SBS^{T}$ and $C=SDS^{T}$;
\emph{unitary congruence}, \emph{*}$\emph{congruence}$, and \emph{unitary
*}$\emph{congruence}$ of two pairs of matrices are defined analogously.

Our consistent point of view is that unitary *congruence is a special kind of
*congruence (rather than a special kind of similarity) that is to be analyzed
with methods from the general theory of *congruence. In a parallel
development, we treat unitary congruence as a special kind of congruence,
rather than as a special kind of consimilarity. \cite[Section 4.6]{HJ1}

The \emph{null space} of a matrix $A$ is denoted by $N(A)=\{x\in\mathbb{C}%
^{n}:Ax=0\}$; $\dim N(A)$, the dimension of $N(A)$, is the \emph{nullity} of
$A$. The quantities $\dim N(A)$, $\dim N(A^{T})$, $\dim\left(  N(A)\cap
N(A^{T})\right)  $, $\dim N(A^{\ast})$, and $\dim\left(  N(A)\cap N(A^{\ast
})\right)  $ play an important role because of their invariance properties:
$\dim N(A)$, $\dim N(A^{T})$, and $\dim\left(  N(A)\cap N(A^{T})\right)  $ are
invariant under congruence; $\dim N(A)$, $\dim N(A^{\ast})$, and $\dim\left(
N(A)\cap N(A^{\ast})\right)  $ are invariant under *congruence.

Suppose $A,U\in M_{n}$ and $U$ is unitary. A computation reveals that if $A$
is conjugate normal (respectively, congruence normal) then $UAU^{T}$ is
conjugate normal (respectively, congruence normal); if $A$ is normal
(respectively, squared normal) then $UAU^{\ast}$ is normal (respectively,
squared normal). Moreover, if $A\in M_{n}$ and $B\in M_{m}$, one verifies that
$A\oplus B$ is, respectively, conjugate normal, congruence normal, normal, or
squared normal if and only if each of $A$ and $B$ has the respective property.

Matrices $A,B$ of the same size (not necessarily square) are \emph{unitarily
equivalent} if there are unitary matrices $V,W$ such that $A=VBW$. Two
matrices are unitarily equivalent if and only if they have the same singular
values, that is, the \emph{singular value decomposition} is a canonical form
for unitary equivalence.

Each $A\in M_{n}$ has a \emph{left} (respectively, \emph{right}) \emph{polar
decomposition} $A=PW$ (respectively, $A$ $=WQ$) in which the Hermitian
positive semidefinite factors $P=(AA^{\ast})^{1/2}$ and $Q=(A^{\ast}A)^{1/2}$
are uniquely determined, $W$ is unitary, and $W=AQ^{-1}=P^{-1}A$ is uniquely
determined if $A$ is nonsingular.

A matrix of the form%
\[
J_{k}(\lambda)=\left[
\begin{array}
[c]{cccc}%
\lambda & 1 &  & 0\\
& \ddots & \ddots & \\
&  & \ddots & 1\\
&  &  & \lambda
\end{array}
\right]  \in M_{k}%
\]
is a \emph{Jordan block with eigenvalue} $\lambda$. The $n$-by-$n$ identity
and zero matrices are denoted by $I_{n}$ and $0_{n}$, respectively.

The \emph{Frobenius norm} of a matrix $A$ is $\left\Vert A\right\Vert
_{F}=\sqrt{\operatorname{tr}\left(  A^{\ast}A\right)  }$: the square root of
the sum of the squares of the absolute values of the entries of $A$. The
\textit{spectral norm} of $A$ is its largest singular value.

In matters of notation and terminology, we follow the conventions in
\cite{HJ1}.

\section{Cosquares, *cosquares, and canonical forms for congruence and
*congruence}

The Jordan Canonical Form of a cosquare or a *cosquare has a very special structure.

\begin{theorem}
[{\cite{HScongruence}, \cite[Theorem 2.3.1]{Wall}}]%
\label{CosquareCharacterize}Let $\mathfrak{A}\in M_{n}$ be nonsingular.
\medskip\medskip%
\newline
(a) $\mathfrak{A}$ is a cosquare if and only if its Jordan Canonical Form is%
\begin{equation}%
{\displaystyle\bigoplus\limits_{k=1}^{\rho}}
\left(  J_{r_{k}}\left(  \left(  -1\right)  ^{r_{k}+1}\right)  \right)  \oplus%
{\displaystyle\bigoplus\limits_{j=1}^{\sigma}}
\left(  J_{s_{j}}\left(  \gamma_{j}\right)  \oplus J_{s_{j}}\left(  \gamma
_{j}^{-1}\right)  \right)  \text{,\quad}\gamma_{j}\in\mathbb{C}\text{, }%
0\neq\gamma_{j}\neq\left(  -1\right)  ^{s_{j}+1}\text{.}\label{JCFcosquare}%
\end{equation}
$\mathfrak{A}$ is a cosquare that is diagonalizable by similarity if and only
if its Jordan Canonical Form is%
\begin{equation}
I\oplus%
{\displaystyle\bigoplus\limits_{j=1}^{q}}
\left[
\begin{array}
[c]{cc}%
\mu_{j}I_{n_{j}} & 0\\
0 & \mu_{j}^{-1}I_{n_{j}}%
\end{array}
\right]  \text{,\quad}\mu_{j}\in\mathbb{C}\text{, }0\neq\mu_{j}\neq
1\text{,}\label{JCFdiagonalizableCosquare}%
\end{equation}
in which $\mu_{1},\mu_{1}^{-1},\ldots,\mu_{q},\mu_{q}^{-1}$ are the distinct
eigenvalues of $\mathfrak{A}$ such that each $\mu_{j}\neq1$; $n_{1}%
,n_{1},\ldots,n_{q},n_{q}$ are their respective multiplicities; the parameters
$\mu_{j}$ in (\ref{JCFdiagonalizableCosquare}) are determined by
$\mathfrak{A}$ up to replacement by $\mu_{j}^{-1}$. \medskip%
\newline
(b) $\mathfrak{A}$ is a *cosquare if and only if its Jordan Canonical Form is%
\begin{equation}
\bigoplus_{k=1}^{\rho}J_{r_{k}}(\beta_{k})\oplus\bigoplus_{j=1}^{\sigma
}\left(  J_{s_{j}}(\gamma_{j})\oplus J_{s_{j}}(\bar{\gamma}_{j}^{-1})\right)
\text{,}\quad\beta_{k},\gamma_{j}\in\mathbb{C}\text{, \ }|\beta_{k}|=1\text{,
}0<\left\vert \gamma_{j}\right\vert <1\text{.}\label{JCF*cosquare}%
\end{equation}
$\mathfrak{A}$ is a *cosquare that is diagonalizable by similarity if and only
if its Jordan Canonical Form is%
\begin{equation}%
{\displaystyle\bigoplus\limits_{k=1}^{p}}
\lambda_{k}I_{m_{k}}\oplus%
{\displaystyle\bigoplus\limits_{j=1}^{q}}
\left[
\begin{array}
[c]{cc}%
\mu_{j}I_{n_{j}} & 0\\
0 & \bar{\mu}_{j}^{-1}I_{n_{j}}%
\end{array}
\right]  \text{,\quad}\lambda_{k},\mu_{j}\in\mathbb{C}\text{, }\left\vert
\lambda_{k}\right\vert =1\text{, }0<\left\vert \mu_{j}\right\vert
<1\text{,}\label{JCFdiagonalizable*Cosquare}%
\end{equation}
in which $\mu_{1},\bar{\mu}_{1}^{-1},\ldots,\mu_{q},\bar{\mu}_{q}^{-1}$ are
the distinct eigenvalues of $\mathfrak{A}$ such that each $\left\vert \mu
_{j}\right\vert \in(0,1)$; $n_{1},n_{1},\ldots,n_{q},n_{q}$ are their
respective multiplicities. The distinct unimodular eigenvalues of
$\mathfrak{A}$ are $\lambda_{1},\ldots,\lambda_{p}$ and their respective
multiplicities are $m_{1},\ldots,m_{p}$.
\end{theorem}

The following theorem involves three types of blocks
\begin{equation}
\Gamma_{k}=%
\begin{bmatrix}
0 &  &  &  & (-1)^{k+1}\\
&  &  & \text{%
\begin{picture}(12,8) \put(-2,-4){$\cdot$} \put(3,0){$\cdot$} \put
(8,4){$\cdot$} \end{picture}%
} & (-1)^{k}\\
&  & 1 & \text{%
\begin{picture}(12,8) \put(-2,-4){$\cdot$} \put(3,0){$\cdot$} \put
(8,4){$\cdot$} \end{picture}%
} & \\
& -1 & -1 &  & \\
1 & 1 &  &  & 0
\end{bmatrix}
\in M_{k},\quad\text{(}\Gamma_{1}=[1]\text{),}\label{Gamman}%
\end{equation}%
\begin{equation}
\Delta_{k}=%
\begin{bmatrix}
0 &  &  & 1\\
&  & \text{%
\begin{picture}(12,8) \put(-2,-4){$\cdot$} \put(3,0){$\cdot$} \put
(8,4){$\cdot$} \end{picture}%
} & i\\
& 1 & \text{%
\begin{picture}(12,8) \put(-2,-4){$\cdot$} \put(3,0){$\cdot$} \put
(8,4){$\cdot$} \end{picture}%
} & \\
1 & i &  & 0
\end{bmatrix}
\in M_{k},\quad\text{(}\Delta_{1}=[1]\text{),}\label{Deltan}%
\end{equation}
and%
\begin{equation}
H_{2k}(\mu)=%
\begin{bmatrix}
0 & I_{k}\\
J_{k}(\mu) & 0
\end{bmatrix}
\in M_{2k},\quad\text{(}H_{2}(\mu)=\left[
\begin{array}
[c]{cc}%
0 & 1\\
\mu & 0
\end{array}
\right]  \text{).}\label{Hn}%
\end{equation}

\begin{theorem}
[\cite{HScongruence}]\label{CongruenceCanonicalForms}Let $A\in M_{n}$ be
nonsingular.
\newline
(a) $A$ is congruent to a direct sum, uniquely determined up to permutation of
summands, of the form%
\begin{equation}%
{\displaystyle\bigoplus\limits_{k=1}^{\rho}}
\Gamma_{r_{k}}\oplus%
{\displaystyle\bigoplus\limits_{j=1}^{\sigma}}
H_{2s_{j}}\left(  \gamma_{j}\right)  ,\quad\gamma_{j}\in\mathbb{C}\text{,
}0\neq\gamma_{j}\neq(-1)^{s_{j}+1}\text{,}\label{ccf}%
\end{equation}
in which each $\gamma_{j}$ is determined up to replacement by $\gamma_{j}%
^{-1}$. If (\ref{JCFcosquare}) is the Jordan Canonical Form of $A^{-T}A$, then
the direct summands in (\ref{ccf}) can be arranged so that the parameters
$\rho$, $\sigma$, $r_{k}$, $s_{j}$, and $\gamma_{j}$ in (\ref{ccf}) are
identical to the same parameters in (\ref{JCFcosquare}). Two nonsingular
matrices are congruent if and only if their cosquares are similar.
\newline
(b) $A$ is *congruent to a direct sum, uniquely determined up to permutation
of summands, of the form%
\begin{equation}%
{\displaystyle\bigoplus\limits_{k=1}^{\rho}}
\alpha_{k}\Delta_{n_{k}}\oplus%
{\displaystyle\bigoplus\limits_{j=1}^{\sigma}}
H_{2m_{j}}\left(  \gamma_{j}\right)  \text{,\quad}\alpha_{k},\gamma_{j}%
\in\mathbb{C}\text{, }|\alpha_{k}|=1\text{, }0<|\gamma_{j}|<1\text{,
}\label{*ccf}%
\end{equation}
If (\ref{JCF*cosquare}) is the Jordan Canonical Form of $A^{-\ast}A$, then the
direct summands in (\ref{*ccf}) can be arranged so that the parameters $r_{k}%
$, $s_{j}$, and $\gamma_{j}$ in (\ref{*ccf}) are identical to the same
parameters in (\ref{JCF*cosquare}), and the parameters $\alpha_{k}$ in
(\ref{*ccf}) and $\beta_{k}$ in (\ref{JCF*cosquare}) satisfy $\alpha_{k}%
^{2}=\beta_{k}$ for each $k=1,\ldots,r$.
\end{theorem}

Among many applications of the canonical form (\ref{*ccf}), it follows that
any complex square matrix is *congruent to its transpose, and the *congruence
can be achieved via a coninvolutory matrix. This conclusion is actually valid
for any square matrix over any field of characteristic not two with an
involution (possibly the identity involution). \cite{HStranspose}

If $A$ is nonsingular and $U$ is unitary, then%
\[
\left(  UAU^{T}\right)  ^{-T}\left(  UAU^{T}\right)  =\bar{U}\left(
A^{-T}A\right)  \bar{U}^{\ast}%
\]
and%
\[
\left(  UAU^{\ast}\right)  ^{-\ast}\left(  UAU^{\ast}\right)  =U\left(
A^{-\ast}A\right)  U^{\ast}\text{,}%
\]
so a unitary congruence (respectively, a unitary *congruence) of a nonsingular
matrix corresponds to a unitary similarity of its cosquare (respectively,
*cosquare), both via the same unitary matrix. If the cosquare or *cosquare of
$A\in M_{n}$ is diagonalizable by unitary similarity, what can be said about a
canonical form for $A$ under unitary congruence or unitary *congruence?{}

\section{Normal matrices, intertwining, and zero blocks}

Intertwining identities involving normal matrices lead to characterizations
and canonical forms for unitary congruences.

\begin{lemma}
\label{Fuglede}Let $A,L,P\in M_{n}$ and assume that $L$ and $P$ are normal.
Then \medskip%
\newline
(a) $AL=PA$ if and only if $AL^{\ast}=P^{\ast}A$. \medskip%
\newline
(b) If $L$ and $P$ are nonsingular, then $AL=PA$ if and only if $AL^{-\ast
}=P^{-\ast}A$.
\end{lemma}

\begin{proof}
Let $L=U\Lambda U^{\ast}$ and $P=V\Pi V^{\ast}$ for some unitary $U,V$ and
diagonal $\Lambda,\Pi$. The intertwining condition $AL=PA$ implies that
$Ag(L)=g(P)A$ for any polynomial $g(t)$. \medskip%
\newline
(a) Let $g(t)$ be any polynomial such that $g(\Lambda)=\bar{\Lambda}$ and
$g(\Pi)=\bar{\Pi}$, that is, $g(t)$ interpolates the function $z\rightarrow
\bar{z}$ on the spectra of $L$ and $P$. Then%
\[
AL^{\ast}=Ag(L)=g(P)A=P^{\ast}A\text{.}%
\]
(b) Use the same argument, but let $g(t)$ interpolate the function
$z\rightarrow\bar{z}^{-1}$ on the spectra of $L$ and $P$.\hfill
\end{proof}

\medskip The following lemma reveals fundamental patterns in the zero blocks
of a partitioned matrix that is normal, conjugate normal, squared normal, or
congruence normal.

\begin{lemma}
\label{Zero Blocks}Let $A\in M_{n}$ be given. \medskip%
\newline
(a) Suppose%
\[
A=\left[
\begin{array}
[c]{cc}%
A_{11} & A_{12}\\
0 & A_{22}%
\end{array}
\right]  \text{,}%
\]
in which $A_{11}$ and $A_{22}$ are square. If $A$ is normal or conjugate
normal, then%
\[
A=\left[
\begin{array}
[c]{cc}%
A_{11} & 0\\
0 & A_{22}%
\end{array}
\right]  \text{.}%
\]
If $A$ is normal, then $A_{11}$ and $A_{22}$ are normal; if $A$ is conjugate
normal, then $A_{11}$ and $A_{22}$ are conjugate normal. \medskip%
\newline
(b) Suppose%
\begin{equation}
A=\left[
\begin{array}
[c]{ccc}%
A_{11} & A_{12} & 0\\
A_{21} & A_{22} & A_{23}\\
0 & 0 & 0_{k}%
\end{array}
\right]  \text{,}\label{generalReduction}%
\end{equation}
in which $A_{11}$ and $A_{22}$ are square, and both $[A_{11}~A_{12}]$ and
$A_{23}$ have full row rank. If $A$ is squared normal or congruence normal,
then%
\[
A=\left[
\begin{array}
[c]{ccc}%
A_{11} & 0 & 0\\
0 & 0 & A_{23}\\
0 & 0 & 0_{k}%
\end{array}
\right]
\]
and $A_{11}$ is nonsingular. If $A$ is squared normal, then $A_{11}$ is
squared normal; if $A$ is congruence normal, then $A_{11}$ is congruence normal.
\end{lemma}

\begin{proof}
(a) If $A$ is normal, then%
\[
A^{\ast}A=\left[
\begin{array}
[c]{cc}%
A_{11}^{\ast}A_{11} & \bigstar\\
\bigstar & \bigstar
\end{array}
\right]  =\left[
\begin{array}
[c]{cc}%
A_{11}A_{11}^{\ast}+A_{12}A_{12}^{\ast} & \bigstar\\
\bigstar & \bigstar
\end{array}
\right]  =AA^{\ast}\text{.}%
\]
We have $A_{11}^{\ast}A_{11}=A_{11}A_{11}^{\ast}+A_{12}A_{12}^{\ast}$, so
$\operatorname{tr}$ $\left(  A_{11}^{\ast}A_{11}\right)  =\operatorname{tr}%
\left(  A_{11}A_{11}^{\ast}\right)  =\operatorname{tr}\left(  A_{11}%
A_{11}^{\ast}\right)  +\operatorname{tr}\left(  A_{12}A_{12}^{\ast}\right)  $.
Then $\operatorname{tr}\left(  A_{12}A_{12}^{\ast}\right)  =\left\Vert
A_{12}^{\ast}\right\Vert _{F}^{2}=0$, so $A_{12}=0$ and $A_{11}^{\ast}%
A_{11}=A_{11}A_{11}^{\ast}$. If $A$ is conjugate normal, then
\[
\overline{A^{\ast}A}=\left[
\begin{array}
[c]{cc}%
\overline{A_{11}^{\ast}A_{11}} & \bigstar\\
\bigstar & \bigstar
\end{array}
\right]  =\left[
\begin{array}
[c]{cc}%
A_{11}A_{11}^{\ast}+A_{12}A_{12}^{\ast} & \bigstar\\
\bigstar & \bigstar
\end{array}
\right]  =AA^{\ast}\text{.}%
\]
We have $\operatorname{tr}\left(  \overline{A_{11}^{\ast}A_{11}}\right)
=\operatorname{tr}\left(  A_{11}A_{11}^{\ast}\right)  =\operatorname{tr}%
\left(  A_{11}A_{11}^{\ast}\right)  +\operatorname{tr}\left(  A_{12}%
A_{12}^{\ast}\right)  $, so $\operatorname{tr}\left(  A_{12}A_{12}^{\ast
}\right)  =\left\Vert A_{12}^{\ast}\right\Vert _{F}^{2}=0$. Then $A_{12}=0$
and $\overline{A_{11}^{\ast}A_{11}}=A_{11}A_{11}^{\ast}$. \medskip%
\newline
(b) Compute
\begin{equation}
A^{2}=\left[
\begin{array}
[c]{ccc}%
\bigstar & \bigstar & A_{12}A_{23}\\
\bigstar & \bigstar & A_{22}A_{23}\\
0 & 0 & 0_{k}%
\end{array}
\right]  \text{ and }\bar{A}A=\left[
\begin{array}
[c]{ccc}%
\bigstar & \bigstar & \overline{A_{12}}A_{23}\\
\bigstar & \bigstar & \overline{A_{22}}A_{23}\\
0 & 0 & 0_{k}%
\end{array}
\right]  \text{.}\label{squares}%
\end{equation}
If $A$ is squared normal (or congruence normal), then (a) ensures that both
$A_{12}A_{23}$ and $A_{22}A_{23}$ (or both $\overline{A_{12}}A_{23}$ and
$\overline{A_{22}}A_{23}$) are zero blocks; since $A_{23}$ has full row rank,
it follows that both $A_{12}$ and $A_{22}$ are zero blocks and hence
\[
A=\left[
\begin{array}
[c]{ccc}%
A_{11} & 0 & 0\\
A_{21} & 0 & A_{23}\\
0 & 0 & 0_{k}%
\end{array}
\right]  \text{,}%
\]
in which $A_{11}$ is nonsingular. Now compute%
\[
A^{2}=\left[
\begin{array}
[c]{cc}%
A_{11}^{2} & 0\\
A_{21}A_{11} & 0
\end{array}
\right]  \oplus0_{k}\text{ and }\bar{A}A=\left[
\begin{array}
[c]{cc}%
\overline{A_{11}}A_{11} & 0\\
\overline{A_{21}}A_{11} & 0
\end{array}
\right]  \oplus0_{k}\text{.}%
\]
If $A$ is squared normal (or if $A$ is congruence normal), then (a) ensures
that $A_{21}A_{11}=0$ (or that $\overline{A_{21}}A_{11}=0$); since $A_{11}$ is
nonsingular, it follows that $A_{21}=0$ and $A_{11}$ is squared normal (or
congruence normal). \hfill
\end{proof}

\medskip A matrix $A\in M_{n}$ is said to be \emph{range Hermitian} if $A$ and
$A^{\ast}$ have the same range. If $\operatorname{rank}A=r$ and there is a
unitary $U$ and a nonsingular $C\in M_{r}$ such that $U^{\ast}AU=C\oplus
0_{n-r}$, then $A$ is range Hermitian; the converse assertion follows from
Theorem \ref{UnitaryRegularization}(b). For example, every normal matrix is
range Hermitian. The following lemma shows that, for a range Hermitian matrix
and a normal matrix, commutativity follows from a generally weaker condition.

\begin{lemma}
\label{Commutivity Implication}Let $A,B\in M_{n}$. Suppose $A$ is range
Hermitian and $B$ is normal. Then $ABA=A^{2}B$ if and only if $AB=BA$.
\end{lemma}

\begin{proof}
If $AB=BA$, then $A(BA)=A(AB)=A^{2}B$. Conversely, suppose $ABA=A^{2}B$. Let
$A=U(C\oplus0_{n-r})U^{\ast}$, in which $U\in M_{n}$ is unitary and $C\in
M_{r}$ is nonsingular. Partition the normal matrix $U^{\ast}BU=[B_{ij}%
]_{i,j=1}^{2}$ conformally to $C\oplus0_{n-r}$. Then%
\begin{align*}
U^{\ast}(ABA)U  & =(U^{\ast}AU)(U^{\ast}BU)(U^{\ast}AU)\\
& =\left[
\begin{array}
[c]{cc}%
C & 0\\
0 & 0
\end{array}
\right]  \left[
\begin{array}
[c]{cc}%
B_{11} & B_{12}\\
B_{21} & B_{22}%
\end{array}
\right]  \left[
\begin{array}
[c]{cc}%
C & 0\\
0 & 0
\end{array}
\right]  =\left[
\begin{array}
[c]{cc}%
CB_{11}C & 0\\
0 & 0
\end{array}
\right]
\end{align*}
and%
\begin{align*}
U^{\ast}(A^{2}B)U  & =(U^{\ast}AU)^{2}(U^{\ast}BU)\\
& =\left[
\begin{array}
[c]{cc}%
C^{2} & 0\\
0 & 0
\end{array}
\right]  \left[
\begin{array}
[c]{cc}%
B_{11} & B_{12}\\
B_{21} & B_{22}%
\end{array}
\right]  =\left[
\begin{array}
[c]{cc}%
C^{2}B_{11} & C^{2}B_{12}\\
0 & 0
\end{array}
\right]  \text{,}%
\end{align*}
so $C^{2}B_{12}=0$, which implies that $B_{12}=0$. Lemma \ref{Zero Blocks}(a)
ensures that $B_{21}=0$ as well, so $U^{\ast}BU=B_{11}\oplus B_{22}$.
Moreover, $CB_{11}C=C^{2}B_{11}$, so $B_{11}C=CB_{11}$. We conclude that
$U^{\ast}AU$ commutes with $U^{\ast}BU$, and hence $A$ commutes with
$B$.\hfill
\end{proof}

\section{Normal cosquares and *cosquares}

A nonsingular matrix $A$ whose cosquare is normal (respectively, whose
*cosquare is normal) has a simple canonical form under unitary congruence
(respectively, under unitary *congruence). Moreover, normality of the cosquare
or *cosquare of $A$ is equivalent to simple properties of $A$ itself that are
the key---via regularization---to obtaining canonical forms under unitary
congruence or unitary *congruence even when $A$ is singular.

\subsection{Normal cosquares}

If $A\in M_{n}$ is nonsingular and its cosquare $\mathfrak{A}$ is normal, then
$\mathfrak{A}$ is unitarily diagonalizable and we may assume that its Jordan
Canonical Form has the form (\ref{JCFdiagonalizableCosquare}). For our
analysis it is convenient to separate the eigenvalue pairs $\{-1,-1\}$ of
$\mathfrak{A}$ from the reciprocal pairs of its other eigenvalues in
(\ref{JCFdiagonalizableCosquare}). Any unitary similarity that puts
$\mathfrak{A}$ in the diagonal form (\ref{JCFdiagonalizableCosquare}) induces
a unitary congruence of $A$ that puts it into a special block diagonal form.

\begin{theorem}
\label{NormalCosquareDiagonalBlocks}Let $A\in M_{n}$ be nonsingular and
suppose that its cosquare $\mathfrak{A}=A^{-T}A$ is normal. Let $\mu_{1}%
,\mu_{1}^{-1},\ldots,\mu_{q},\mu_{q}^{-1}$ be the distinct eigenvalues of
$\mathfrak{A}$ with $-1\neq\mu_{j}\neq1$ for each $j=1,\ldots,q$, and let
$n_{1},n_{1},\ldots,n_{q},n_{q}$ be their respective multiplicities. Let
$n_{+}$ and $2n_{-}$be the multiplicities of $+1$ and $-1$, respectively, as
eigenvalues of $\mathfrak{A}$. Let
\begin{equation}
\Lambda=I_{n_{+}}\oplus\left(  -I_{2n_{-}}\right)  \oplus%
{\displaystyle\bigoplus\limits_{j=1}^{q}}
\left[
\begin{array}
[c]{cc}%
\mu_{j}I_{n_{j}} & 0\\
0 & \mu_{j}^{-1}I_{n_{j}}%
\end{array}
\right]  \text{,\quad}\mu_{j}\neq0\text{, }-1\neq\mu_{j}\neq1\text{,}%
\label{NormalCosquareLambda}%
\end{equation}
let $U\in M_{n}$ be any unitary matrix such that $\mathfrak{A}=U\Lambda
U^{\ast}$, and let $\mathcal{A}=U^{T}AU$. Then%
\begin{equation}
\mathcal{A}=\mathcal{A}_{+}\oplus\mathcal{A}_{-}\oplus\mathcal{A}_{1}%
\oplus\cdots\oplus\mathcal{A}_{q}\text{,}\label{NormalCosquareBlocks}%
\end{equation}
in which $\mathcal{A}_{+}\in M_{n_{+}}$ is symmetric, $\mathcal{A}_{-}\in
M_{2n_{-}}$ is skew symmetric, and each $\mathcal{A}_{j}\in M_{2n_{j}}$ has
the form%
\begin{equation}
\mathcal{A}_{j}=\left[
\begin{array}
[c]{cc}%
0_{n_{j}} & Y_{j}\\
\mu_{j}Y_{j}^{T} & 0_{n_{j}}%
\end{array}
\right]  \text{,\quad}Y_{j}\in M_{n_{j}}\text{ is nonsingular.}%
\label{NormalCosquareThirdTypeBlock}%
\end{equation}
The unitary congruence class of each of the $q+2$ blocks in
(\ref{NormalCosquareBlocks}) is uniquely determined.
\end{theorem}

\begin{proof}
The presentation (\ref{NormalCosquareLambda}) of the Jordan Canonical Form of
$\mathfrak{A}$ differs from that in (\ref{JCFdiagonalizableCosquare}) only in
the separate identification of the eigenvalue pairs $\{-1,-1\}$. We have
$A=A^{T}\mathfrak{A}=A^{T}U\Lambda U^{\ast}$, which implies that%
\[
\mathcal{A}=U^{T}AU=U^{T}A^{T}U\Lambda=\mathcal{A}^{T}\Lambda
\]
and hence%
\[
\mathcal{A}=\mathcal{A}^{T}\Lambda=\left(  \mathcal{A}^{T}\Lambda\right)
^{T}\Lambda=\Lambda\mathcal{A}\Lambda\text{,}%
\]
that is,%
\begin{equation}
\Lambda^{-1}\mathcal{A}=\mathcal{A}\Lambda\text{.}\label{CosquareCommute}%
\end{equation}
Partition $\mathcal{A}=\left[  \mathcal{A}_{ij}\right]  _{i,j=1}^{q+2}$
conformally to $\Lambda$. The $q+2$ diagonal blocks of $\Lambda$ have mutually
distinct spectra; the spectra of corresponding diagonal blocks of $\Lambda$
and $\Lambda^{-1}$ are the same. The identity (\ref{CosquareCommute}) and
Sylvester's Theorem on Linear Matrix Equations \cite[Section 2.4, Problems 9
and 13]{HJ1} ensure that $\mathcal{A}$ is block diagonal and conformal to
$\Lambda$, that is,
\[
\mathcal{A}=\mathcal{A}_{11}\oplus\cdots\oplus\mathcal{A}_{q+2,q+2}%
\]
is block diagonal. Moreover, the identity $\mathcal{A}=\mathcal{A}^{T}\Lambda$
ensures that (a) $\mathcal{A}_{11}=\mathcal{A}_{11}^{T}$, so $\mathcal{A}%
_{+}:=\mathcal{A}_{11}$ is symmetric; (b) $\mathcal{A}_{22}=-\mathcal{A}%
_{22}^{T}$, so $\mathcal{A}_{-}:=\mathcal{A}_{22}$ is skew symmetric; and (c)
for each $i=3,\ldots,q+2$ the nonsingular block $\mathcal{A}_{jj}$ has the
form%
\[
\left[
\begin{array}
[c]{cc}%
X & Y\\
Z & W
\end{array}
\right]  \text{,\quad}X,Y,Z,W\in M_{n_{j}}%
\]
and satisfies an identity of the form%
\[
\left[
\begin{array}
[c]{cc}%
X & Y\\
Z & W
\end{array}
\right]  =\left[
\begin{array}
[c]{cc}%
X & Y\\
Z & W
\end{array}
\right]  ^{T}\left[
\begin{array}
[c]{cc}%
\mu I & 0\\
0 & \mu^{-1}I
\end{array}
\right]  =\left[
\begin{array}
[c]{cc}%
\mu X^{T} & \mu^{-1}Z^{T}\\
\mu Y^{T} & \mu^{-1}W^{T}%
\end{array}
\right]
\]
in which $\mu^{2}\neq1$. But $X=\mu X^{T}=\mu^{2}X$ and $W=\mu^{-1}W^{T}%
=\mu^{-2}W$, so $X=W=0$. Moreover, $Z=\mu Y^{T}$, so $\mathcal{A}_{jj}$ has
the form (\ref{NormalCosquareThirdTypeBlock}).

What can we say if $\mathfrak{A}$ can be put into the form
(\ref{NormalCosquareLambda}) via unitary similarity with a different unitary
matrix $V$? If $\mathfrak{A}=U\Lambda U^{\ast}=V\Lambda V^{\ast}$ and both $U$
and $V$ are unitary, then $\Lambda\left(  U^{\ast}V\right)  =\left(  U^{\ast
}V\right)  \Lambda$, so another application of Sylvester's Theorem ensures
that the unitary matrix $U^{\ast}V$ is block diagonal and conformal to
$\Lambda$. Thus, in the respective presentations (\ref{NormalCosquareBlocks})
associated with $U$ and $V$, corresponding diagonal blocks are unitarily
congruent.\hfill
\end{proof}

\begin{theorem}
\label{NonsingularEquivalence}Let $A\in M_{n}$. The following are equivalent:
\newline
(a) $\bar{A}A$ is normal.
\newline
(b) $A\left(  \overline{AA^{\ast}}\right)  =\left(  \overline{A^{\ast}%
A}\right)  A$, that is, $A\bar{A}A^{T}=A^{T}\bar{A}A$. \medskip%
\newline
If $A$ is nonsingular, then (a) and (b) are equivalent to%
\newline
(c) $A^{-T}A$ is normal.
\end{theorem}

\begin{proof}
(a) $\Rightarrow$ (b): Consider the identity%
\[
A\left(  \bar{A}A\right)  =A\bar{A}A=\left(  A\bar{A}\right)  A\text{.}%
\]
Since $\bar{A}A$ is normal, $A\bar{A}=\overline{(\bar{A}A)}$ is normal and
Lemma \ref{Fuglede}(a) ensures that%
\[
A\left(  \bar{A}A\right)  ^{\ast}=AA^{\ast}A^{T}=A^{T}A^{\ast}A=\left(
A\bar{A}\right)  ^{\ast}A\text{.}%
\]
Taking the transpose of the middle identity, and using Hermicity of $AA^{\ast
}$ and $A^{\ast}A$, gives%
\[
A\left(  AA^{\ast}\right)  ^{T}=A\left(  \overline{AA^{\ast}}\right)  =\left(
\overline{A^{\ast}A}\right)  A=\left(  A^{\ast}A\right)  ^{T}A\text{.}%
\]
(b) $\Rightarrow$ (a): We use (b) in the form $AA^{\ast}A^{T}=A^{T}A^{\ast}A$
to compute
\[
\left(  \bar{A}A\right)  \left(  \bar{A}A\right)  ^{\ast}=\bar{A}\left(
AA^{\ast}A^{T}\right)  =\bar{A}\left(  A^{T}A^{\ast}A\right)  \text{,}%
\]
so $\bar{A}A^{T}A^{\ast}A$ is Hermitian:%
\begin{align*}
\left(  \bar{A}A\right)  \left(  \bar{A}A\right)  ^{\ast}  & =\bar{A}%
A^{T}A^{\ast}A=\left(  \bar{A}A^{T}A^{\ast}A\right)  ^{\ast}\\
& =A^{\ast}\left(  A\bar{A}A^{T}\right)  =A^{\ast}\left(  A^{T}\bar
{A}A\right)  =\left(  \bar{A}A\right)  ^{\ast}\left(  \bar{A}A\right)
\text{.}%
\end{align*}
(c) $\Rightarrow$ (b): Consider the identity%
\[
A(A^{-T}A)=AA^{-T}A=(A^{-T}A)^{-T}A\text{.}%
\]
Since $A^{-T}A$ is normal, Lemma \ref{Fuglede}(b) ensures that%
\[
A(A^{-T}A)^{-\ast}=\left(  (A^{-T}A)^{-T}\right)  ^{-\ast}A=\overline
{(A^{-T}A)}A\text{,}%
\]
so $A\bar{A}A^{-\ast}=A^{-\ast}\bar{A}A$, from which it follows that $A^{\ast
}A\bar{A}=\bar{A}AA^{\ast}$, which is the conjugate of (b).\medskip%
\newline
(b) $\Rightarrow$ (c): Since $A$ is nonsingular, the identity (b) is
equivalent to%
\[
A^{T}\bar{A}A=\left(  A^{T}A^{\ast}A\right)  ^{T}=\left(  AA^{\ast}%
A^{T}\right)  ^{T}=A\bar{A}A^{T}\text{,}%
\]
which in turn is equivalent to%
\[
AA^{-T}\bar{A}^{-1}=\bar{A}^{-1}A^{-T}A\text{.}%
\]
Now compute%
\begin{align*}
(A^{-T}A)(A^{-T}A)^{\ast}  & =A^{-T}\left(  AA^{\ast}\right)  \bar{A}%
^{-1}=A^{-T}\left(  A^{T}A^{\ast}AA^{-T}\right)  \bar{A}^{-1}\\
& =A^{\ast}\left(  AA^{-T}\bar{A}^{-1}\right)  =A^{\ast}\left(  \bar{A}%
^{-1}A^{-T}A\right)  =(A^{-T}A)^{\ast}(A^{-T}A)\text{.}%
\end{align*}
\hfill
\end{proof}

\begin{theorem}
\label{nonsingularTheorem}Let $A\in M_{n}$ be nonsingular. If $\bar{A}A$ is
normal, then $A$ is unitarily congruent to a direct sum of blocks, each of
which is%
\begin{equation}
\left[  \sigma\right]  \text{ or }\tau\left[
\begin{array}
[c]{cc}%
0 & 1\\
\mu & 0
\end{array}
\right]  \text{,\quad}\sigma>0\text{, }\tau>0\text{, }\mu\in\mathbb{C}\text{,
}0\neq\mu\neq1\text{.}\label{cnc0}%
\end{equation}
This direct sum is uniquely determined by $A$ up to permutation of its blocks
and replacement of any $\mu$ by $\mu^{-1}$. Conversely, if $A$ is unitarily
congruent to a direct sum of blocks of the two types in (\ref{cnc0}), then
$\bar{A}A$ is normal.
\end{theorem}

\begin{proof}
Normality of $\bar{A}A$ implies normality of the cosquare $A^{-T}A$. Theorem
\ref{NormalCosquareDiagonalBlocks} ensures that $A$ is unitarily congruent to
a direct sum of the form (\ref{NormalCosquareBlocks}), and the unitary
congruence class of each summand is uniquely determined by $A$. It suffices to
consider the three types of blocks that occur in (\ref{NormalCosquareBlocks}):
(a) a symmetric block $\mathcal{A}_{+}$, (b) a skew-symmetric block
$\mathcal{A}_{-}$, and (c) a block of the form
(\ref{NormalCosquareThirdTypeBlock}). \medskip%
\newline
(a) The special singular value decomposition available for a nonsingular
symmetric matrix \cite[Corollary 4.4.4]{HJ1} ensures that there is a unitary
$V$ and a positive diagonal matrix $\Sigma=\operatorname{diag}(\sigma
_{1},\ldots,\sigma_{n_{+}})$ such that $\mathcal{A}_{+}=V\Sigma V^{T}$. The
singular values $\sigma_{i}$ of $\mathcal{A}_{+}$ are the source of all of the
1-by-1 blocks in (\ref{cnc0}). They are unitary congruence invariants of
$\mathcal{A}_{+}$, so they are uniquely determined by $A$. \medskip%
\newline
(b) The special singular value decomposition available for a nonsingular
skew-symmetric matrix \cite[Problem 26, Section 4.4]{HJ1} ensures that there
is a unitary $V$ and a nonsingular block diagonal matrix%
\begin{equation}
\Sigma=\tau_{1}\left[
\begin{array}
[c]{cc}%
0 & 1\\
-1 & 0
\end{array}
\right]  \oplus\cdots\oplus\tau_{n_{-}}\left[
\begin{array}
[c]{cc}%
0 & 1\\
-1 & 0
\end{array}
\right] \label{skewsymmetric}%
\end{equation}
such that $\mathcal{A}_{-}=V\Sigma V^{T}$. These blocks are the source of all
of the 2-by-2 blocks in (\ref{cnc0}) in which $\mu=-1$. The parameters
$\tau_{1},\tau_{1},\ldots,\tau_{n_{-}},\tau_{n_{-}}$ are the singular values
of $\mathcal{A}_{-}$, which are unitary congruence invariants of
$\mathcal{A}_{-}$, so they are uniquely determined by $A$. \medskip%
\newline
(c) Consider a block of the form%
\[
\mathcal{A}_{j}=\left[
\begin{array}
[c]{cc}%
0 & Y_{j}\\
\mu_{j}Y_{j}^{T} & 0
\end{array}
\right]
\]
in which $Y_{j}\in M_{n_{j}}$ is nonsingular. The singular value decomposition
\cite[Theorem 7.3.5]{HJ1} ensures that there are unitary $V_{j},W_{j}\in
M_{n_{j}}$ and a positive diagonal matrix $\Sigma_{j}=\operatorname{diag}%
(\tau_{1}^{(j)},\ldots,\tau_{n_{j}}^{(j)})$ such that $Y_{j}=V_{j}\Sigma
_{j}W_{j}^{\ast}$. Then%
\[
\mathcal{A}_{j}=\left[
\begin{array}
[c]{cc}%
0 & V_{j}\Sigma_{j}W_{j}^{\ast}\\
\mu_{j}\bar{W}_{j}\Sigma_{j}V_{j}^{T} & 0
\end{array}
\right]  =\left[
\begin{array}
[c]{cc}%
V_{j} & 0\\
0 & \bar{W}_{j}%
\end{array}
\right]  \left[
\begin{array}
[c]{cc}%
0 & \Sigma_{j}\\
\mu_{j}\Sigma_{j} & 0
\end{array}
\right]  \left[
\begin{array}
[c]{cc}%
V_{j} & 0\\
0 & \bar{W}_{j}%
\end{array}
\right]  ^{T}%
\]
is unitarily congruent to%
\[
\left[
\begin{array}
[c]{cc}%
0 & \Sigma_{j}\\
\mu_{j}\Sigma_{j} & 0
\end{array}
\right]  \text{,}%
\]
which is unitarily congruent (permutation similar) to%
\[%
{\displaystyle\bigoplus\limits_{i=1}^{n_{j}}}
\tau_{i}^{(j)}\left[
\begin{array}
[c]{cc}%
0 & 1\\
\mu_{j} & 0
\end{array}
\right]  \text{,\quad}\tau_{i}^{(j)}>0\text{.}%
\]
These blocks contribute $n_{j}$ 2-by-2 blocks to (\ref{cnc0}), all with
$\mu=\mu_{j}$. Given $\mu_{j}\neq0$, the parameters $\tau_{1}^{(j)}%
,\ldots,\tau_{n_{j}}^{(j)}$ are determined by the eigenvalues of
$\overline{\mathcal{A}_{j}}\mathcal{A}_{j}$, which are invariant under unitary
congruence of $\mathcal{A}_{j}$.

Conversely, if $A$ is unitarily congruent to a direct sum of blocks of the
form (\ref{cnc0}), then $\bar{A}A$ is unitarily similar to a direct sum of
blocks, each of which is%
\[
\lbrack\sigma^{2}]\text{ or }\tau^{2}\mu I_{2}\text{,\quad\ }0\neq\mu
\neq1\text{,}%
\]
so $\bar{A}A$ is normal.\hfill
\end{proof}

\subsection{Normal *cosquares}

If $A\in M_{n}$ is nonsingular and its *cosquare $\mathfrak{A}$ is normal, we
can deduce a unitary *congruence canonical form for $A$ with an argument
largely parallel to that in the preceding section, starting with the Jordan
Canonical Form (\ref{JCFdiagonalizable*Cosquare}). We find that any unitary
similarity that diagonalizes $\mathfrak{A}$ induces a unitary *congruence of
$A$ that puts it into a special block diagonal form.

\begin{theorem}
\label{Normal*CosquareDiagonalBlocks}Let $A\in M_{n}$ be nonsingular and
suppose that its *cosquare $\mathfrak{A}=A^{-\ast}A$ is normal. Let $\mu
_{1},\bar{\mu}_{1}^{-1},\ldots,\mu_{q},\bar{\mu}_{q}^{-1}$ be the distinct
eigenvalues of $\mathfrak{A}$ with $0<|\mu_{j}|<1$ for each $j=1,\ldots,q$,
and let $n_{1},n_{1},\ldots,n_{q},n_{q}$ be their respective multiplicities.
Let $\lambda_{1},\ldots,\lambda_{p}$ be the distinct unimodular eigenvalues of
$\mathfrak{A}$, with respective multiplicities $m_{1},\ldots,m_{p}$, and
choose any unimodular parameters $\alpha_{1},\ldots,\alpha_{p}$ such that
$\alpha_{k}^{2}=\lambda_{k}$ for each $k=1,\ldots,p$. Let
\begin{equation}
\Lambda=%
{\displaystyle\bigoplus\limits_{k=1}^{p}}
\lambda_{k}I_{m_{k}}\oplus%
{\displaystyle\bigoplus\limits_{j=1}^{q}}
\left[
\begin{array}
[c]{cc}%
\mu_{j}I_{n_{j}} & 0\\
0 & \overline{\mu_{j}}^{-1}I_{n_{j}}%
\end{array}
\right]  \text{,\quad}\lambda_{k},\mu_{j}\in\mathbb{C}\text{, }\left\vert
\lambda_{k}\right\vert =1\text{, }0<\left\vert \mu_{j}\right\vert
<1\text{,}\label{Normal*CosquareLambda}%
\end{equation}
let $U\in M_{n}$ be any unitary matrix such that $\mathfrak{A}=U\Lambda
U^{\ast}$, and let $\mathcal{A}=U^{\ast}AU$. Then $\mathcal{A}$ is block
diagonal and has the form%
\begin{equation}
\mathcal{A}=\alpha_{1}\mathcal{H}_{1}\oplus\cdots\oplus\alpha_{p}%
\mathcal{H}_{p}\oplus\mathcal{A}_{1}\oplus\cdots\oplus\mathcal{A}_{q}%
\text{,}\label{Normal*CosquareBlocks}%
\end{equation}
in which $\mathcal{H}_{k}\in M_{m_{k}}$ is Hermitian for each $k=1,\ldots,p$,
and each $\mathcal{A}_{j}\in M_{2n_{j}}$ has the form%
\begin{equation}
\mathcal{A}_{j}=\left[
\begin{array}
[c]{cc}%
0 & Y_{j}\\
\mu_{j}Y_{j}^{\ast} & 0
\end{array}
\right]  \text{,\quad}Y_{j}\in M_{n_{j}}\text{ is nonsingular.}%
\label{Normal*CosquareThirdTypeBlock}%
\end{equation}
For a given ordering of the blocks in (\ref{Normal*CosquareLambda}) the
unitary *congruence class of each of the $p+q$ blocks in
(\ref{Normal*CosquareBlocks}) is uniquely determined.
\end{theorem}

\begin{proof}
We have $A=A^{\ast}\mathfrak{A}=A^{\ast}U\Lambda U^{\ast}$, which implies that%
\[
\mathcal{A}=U^{\ast}AU=U^{\ast}A^{\ast}U\Lambda=\mathcal{A}^{\ast}\Lambda
\]
and hence%
\[
\mathcal{A}=\mathcal{A}^{\ast}\Lambda=\left(  \mathcal{A}^{\ast}%
\Lambda\right)  ^{\ast}\Lambda=\bar{\Lambda}\mathcal{A}\Lambda\text{,}%
\]
that is,
\begin{equation}
\bar{\Lambda}^{-1}\mathcal{A}=\mathcal{A}\Lambda\text{.}%
\label{*CosquareCommute}%
\end{equation}
Partition $\mathcal{A}=\left[  \mathcal{A}_{ij}\right]  _{i,j=1}^{p+q}$
conformally to $\Lambda$. The $p+q$ diagonal blocks of $\Lambda$ have mutually
distinct spectra; the spectra of corresponding blocks of $\Lambda$ and
$\bar{\Lambda}^{-1}$ are the same. The identity (\ref{*CosquareCommute}) and
Sylvester's Theorem on Linear Matrix Equations ensure that $\mathcal{A}$ is
block diagonal and conformal to $\Lambda$, that is,
\[
\mathcal{A}=\mathcal{A}_{11}\oplus\cdots\oplus\mathcal{A}_{pp}\oplus
\mathcal{A}_{p+1,p+1}\oplus\cdots\oplus\mathcal{A}_{p+q,p+q}%
\]
is block diagonal. Moreover, the identity $\mathcal{A}=\mathcal{A}^{\ast
}\Lambda$ ensures that $\mathcal{A}_{kk}=\lambda_{k}\mathcal{A}_{kk}^{\ast
}=\alpha_{k}^{2}\mathcal{A}_{kk}^{\ast}$, so if we define $\mathcal{H}%
_{k}:=\overline{\alpha_{k}}\mathcal{A}_{kk}$, then
\[
\mathcal{H}_{k}=\overline{\alpha_{k}}\mathcal{A}_{kk}=\overline{\alpha_{k}%
}\alpha_{k}^{2}\mathcal{A}_{kk}^{\ast}=\alpha_{k}\mathcal{A}_{kk}^{\ast
}=\left(  \overline{\alpha_{k}}\mathcal{A}_{kk}\right)  ^{\ast}=\mathcal{H}%
_{k}^{\ast}\text{,}%
\]
so $\mathcal{H}_{k}$ is Hermitian. For each $j=p+1,\ldots,p+q$ the block
$\mathcal{A}_{jj}$ has the form%
\[
\left[
\begin{array}
[c]{cc}%
X & Y\\
Z & W
\end{array}
\right]
\]
and satisfies an identity of the form%
\[
\left[
\begin{array}
[c]{cc}%
X & Y\\
Z & W
\end{array}
\right]  =\left[
\begin{array}
[c]{cc}%
X & Y\\
Z & W
\end{array}
\right]  ^{\ast}\left[
\begin{array}
[c]{cc}%
\mu I & 0\\
0 & \bar{\mu}^{-1}I
\end{array}
\right]  =\left[
\begin{array}
[c]{cc}%
\mu X^{\ast} & \bar{\mu}^{-1}Z^{\ast}\\
\mu Y^{\ast} & \bar{\mu}^{-1}W^{\ast}%
\end{array}
\right]
\]
in which $\left\vert \mu\right\vert ^{2}>1$. But $X=\mu X^{\ast}=\left\vert
\mu\right\vert ^{2}X$ and $W=\overline{\mu^{-1}}W^{\ast}=\left\vert
\mu\right\vert ^{-2}W$, so $X=W=0$. Moreover, $Z=\mu Y^{\ast}$, so
$\mathcal{A}_{jj}$ has the form (\ref{Normal*CosquareThirdTypeBlock}).

If $\mathfrak{A}=U\Lambda U^{\ast}=V\Lambda V^{\ast}$ and both $U$ and $V$ are
unitary, then $\Lambda\left(  U^{\ast}V\right)  =\left(  U^{\ast}V\right)
\Lambda$, so the unitary matrix $U^{\ast}V$ is block diagonal and conformal to
$\Lambda$. Thus, in the presentations (\ref{Normal*CosquareBlocks})
corresponding to $U$ and to $V$, corresponding diagonal blocks are unitarily
*congruent.\hfill
\end{proof}

\begin{theorem}
\label{Nonsingular*equivalence}Let $A\in M_{n}$. The following are
equivalent:
\newline
(a) $A^{2}$ is normal.
\newline
(b) $A\left(  AA^{\ast}\right)  =\left(  A^{\ast}A\right)  A$, that is,
$A^{2}A^{\ast}=A^{\ast}A^{2}$. \medskip%
\newline
If $A$ is nonsingular, then (a) and (b) are equivalent to
\newline
(c) $A^{-\ast}A$ is normal.
\end{theorem}

\begin{proof}
(a) $\Rightarrow$ (b): Consider the identity%
\[
\left(  A^{2}\right)  A=A\left(  A^{2}\right)  \text{.}%
\]
Since $A^{2}$ is normal, Lemma \ref{Fuglede}(a) ensures that%
\[
\left(  A^{2}\right)  ^{\ast}A=A\left(  A^{2}\right)  ^{\ast}\text{,}%
\]
and hence
\[
A^{2}A^{\ast}=A\left(  AA^{\ast}\right)  =\left(  A^{\ast}A\right)  A=A^{\ast
}A^{2}.
\]
(b) $\Rightarrow$ (a): Assuming (b), we have
\begin{align*}
A^{2}\left(  A^{2}\right)  ^{\ast}  & =\left(  A^{2}A^{\ast}\right)  A^{\ast
}=\left(  A^{\ast}A^{2}\right)  A^{\ast}=A^{\ast}\left(  A^{2}A^{\ast}\right)
\\
& =A^{\ast}\left(  A^{\ast}A^{2}\right)  =\left(  A^{\ast}\right)  ^{2}%
A^{2}=\left(  A^{2}\right)  ^{\ast}A^{2}\text{.}%
\end{align*}
(c) $\Rightarrow$ (b): The identity $A=(A^{-\ast}A)^{\ast}A(A^{-\ast}A)$ is
equivalent to%
\begin{equation}
(A^{-\ast}A)^{-\ast}A=A(A^{-\ast}A)\text{.}\label{squared1}%
\end{equation}
Since $AA^{-\ast}$ is normal, Lemma \ref{Fuglede}(b) ensures that%
\[
(A^{-\ast}A)A=\left(  (A^{-\ast}A)^{-\ast}\right)  ^{-\ast}A=A(A^{-\ast
}A)^{-\ast}\text{,}%
\]
which implies that $A^{-\ast}A^{2}=A^{2}A^{-\ast}$ and $A^{2}A^{\ast}=A^{\ast
}A^{2}$.\medskip%
\newline
(b) $\Rightarrow$ (c): Assuming (b), we have $AAA^{\ast}=A^{\ast}AA$, which is
equivalent to $AA^{\ast}A^{\ast}=A^{\ast}A^{\ast}A$ and (since $A$ is
nonsingular) to%
\[
A^{-\ast}AA^{\ast}=A^{\ast}AA^{-\ast}\text{,}%
\]
The inverse of this identity is%
\[
A^{-\ast}A^{-1}A^{\ast}=A^{\ast}A^{-1}A^{-\ast}\text{.}%
\]
Now compute%
\[
(A^{-\ast}A)(A^{-\ast}A)^{\ast}=\left(  A^{-\ast}AA^{\ast}\right)
A^{-1}=\left(  A^{\ast}AA^{-\ast}\right)  A^{-1}\text{,}%
\]
which is Hermitian, so%
\begin{align*}
(A^{-\ast}A)(A^{-\ast}A)^{\ast}  & =\left(  A^{\ast}AA^{-\ast}A^{-1}\right)
^{\ast}=\left(  A^{-\ast}A^{-1}A^{\ast}\right)  A\\
& =\left(  A^{\ast}A^{-1}A^{-\ast}\right)  A=(A^{-\ast}A)^{\ast}(A^{-\ast
}A)\text{.}%
\end{align*}

\hfill
\end{proof}

\begin{theorem}
\label{*nonsingularTheorem}Let $A\in M_{n}$ be nonsingular. If $A^{2}$ is
normal, then $A$ is unitarily *congruent to a direct sum of blocks, each of
which is%
\begin{equation}
\lbrack\lambda]\text{ or }\tau\left[
\begin{array}
[c]{cc}%
0 & 1\\
\mu & 0
\end{array}
\right]  \text{,\quad}\lambda,\mu\in\mathbb{C}\text{, }\lambda\neq0\text{,
}\tau>0\text{, }0<|\mu|<1\text{.}\label{*nonsingular}%
\end{equation}
This direct sum is uniquely determined by $A$, up to permutation of its
blocks. Conversely, if $A$ is unitarily *congruent to a direct sum of blocks
of the form (\ref{*nonsingular}), then $A^{2}$ is normal.
\end{theorem}

\begin{proof}
Normality of $A^{2}$ implies normality of its *cosquare $A^{-\ast}A$, so $A$
is unitarily *congruent to a direct sum of the form
(\ref{Normal*CosquareBlocks}), and the unitary *congruence class of each
summand is uniquely determined by $A$. It suffices to consider the two types
of blocks that occur in (\ref{Normal*CosquareBlocks}): (a) a unimodular scalar
multiple of a Hermitian matrix $\mathcal{H}_{k}$, and (b) a block of the form
(\ref{Normal*CosquareThirdTypeBlock}). \medskip%
\newline
(a) The spectral theorem ensures that there is a unitary $V_{k}\in M_{m_{k}}$
and a real nonsingular diagonal $L_{k}\in M_{m_{k}}$ such that $\mathcal{H}%
_{k}=V_{k}L_{k}V_{k}^{\ast}$. The diagonal entries of $\alpha_{k}L_{k}$ (that
is, the eigenvalues of $\alpha_{k}\mathcal{H}_{k}=$ $\mathcal{A}_{k}$) are
unitary *congruence invariants of $\mathcal{A}_{k}$, so they are uniquely
determined by $A$; they all lie on the line $\{t\alpha_{k}:-\infty<t<\infty
\}$. The diagonal entries of $\alpha_{1}L_{1},\ldots,\alpha_{p}L_{p}$ are the
source of all the 1-by-1 blocks in (\ref{*nonsingular}). \medskip%
\newline
(b) Consider a block of the form%
\[
\mathcal{A}_{j}=\left[
\begin{array}
[c]{cc}%
0 & Y_{j}\\
\mu_{j}Y_{j}^{\ast} & 0
\end{array}
\right]  \text{,\quad}0<\left\vert \mu_{j}\right\vert <1\text{,}%
\]
in which $Y_{j}\in M_{n_{j}}$ is nonsingular. The singular value decomposition
ensures that there are unitary $V_{j},W_{j}\in M_{n_{j}}$ and a positive
diagonal matrix $\Sigma_{j}=\operatorname{diag}(\tau_{1}^{(j)},\ldots
,\tau_{n_{j}}^{(j)})$ such that $Y_{j}=V_{j}\Sigma_{j}W_{j}^{\ast}$. Then%
\[
\mathcal{A}_{j}=\left[
\begin{array}
[c]{cc}%
0 & V_{j}\Sigma_{j}W_{j}^{\ast}\\
\mu_{j}W_{j}\Sigma_{j}V_{j}^{\ast} & 0
\end{array}
\right]  =\left[
\begin{array}
[c]{cc}%
V_{j} & 0\\
0 & W_{j}%
\end{array}
\right]  \left[
\begin{array}
[c]{cc}%
0 & \Sigma_{j}\\
\mu_{j}\Sigma_{j} & 0
\end{array}
\right]  \left[
\begin{array}
[c]{cc}%
V_{j} & 0\\
0 & W_{j}%
\end{array}
\right]  ^{\ast}%
\]
is unitarily *congruent to%
\[
\left[
\begin{array}
[c]{cc}%
0 & \Sigma_{j}\\
\mu_{j}\Sigma_{j} & 0
\end{array}
\right]  \text{,}%
\]
which is unitarily *congruent (permutation similar) to%
\[%
{\displaystyle\bigoplus\limits_{i=1}^{n_{j}}}
\tau_{i}^{(j)}\left[
\begin{array}
[c]{cc}%
0 & 1\\
\mu_{j} & 0
\end{array}
\right]  \text{,\quad}\tau_{i}^{(j)}>0\text{.}%
\]
These blocks contribute $n_{j}$ 2-by-2 blocks to (\ref{*nonsingular}), all
with $\mu=\mu_{j}$. Given $\mu_{j}\neq0$, the parameters $\tau_{1}%
^{(j)},\ldots,\tau_{n_{j}}^{(j)}$ are determined by the eigenvalues of
$\mathcal{A}_{j}^{2}$, which are invariant under unitary *congruence of
$\mathcal{A}_{j}$.

Conversely, if $A$ is unitarily *congruent to a direct sum of blocks of the
two types in (\ref{*nonsingular}), then $A^{2}$ is normal since it is
unitarily *congruent to a direct sum of diagonal blocks of the two types
$[\lambda^{2}]$ and $\tau^{2}\mu I_{2}$.\hfill
\end{proof}

\section{Unitary regularization}

The following theorem describes a reduced form that can be achieved for any
singular nonzero matrix under both unitary congruence and unitary *congruence.
It is the key to separating a singular nonzero matrix into a canonical direct
sum of its regular and singular parts under unitary congruence or unitary *congruence.

\begin{theorem}
\label{UnitaryRegularization}Let $A\in M_{n}$ be singular and nonzero, let
$m_{1}$ be the nullity of $A$, let the columns of $V_{1}$ be any orthonormal
basis for the range of $A$, let the columns of $V_{2}$ be any orthonormal
basis for the null space of $A^{\ast}$, and form the unitary matrix
$V=[V_{1}~V_{2}]$. Then \medskip%
\newline
(a) $A$ is unitarily congruent to a reduced form%
\begin{equation}
\left[
\begin{array}
[c]{cc|c}%
A^{\prime} & B & 0\\
C & D & \left[  \Sigma~0\right] \\\hline
\multicolumn{2}{c|}{0} & 0_{m_{1}}%
\end{array}
\right]  \hspace{-0.13in}%
\begin{array}
[c]{l}%
\\
\}m_{2}\\
\}m_{1}%
\end{array}
\label{Reduced Form}%
\end{equation}
in which $m_{2}=m_{1}-\dim\left(  N(A)\cap N(A^{T})\right)  $; if $m_{2}>0$
then $D\in M_{m_{2}}$, $\Sigma=\operatorname{diag}(\sigma_{1},\ldots
,\sigma_{m_{2}})$, and all $\sigma_{i}>0$; if $m_{1}+m_{2}<n$, then
$A^{\prime}$ is square and $[A^{\prime}~B]$ has linearly independent rows; the
integers $m_{1}$, $m_{2}$ and the unitary congruence class of the block%
\begin{equation}
\left[
\begin{array}
[c]{cc}%
A^{\prime} & B\\
C & D
\end{array}
\right] \label{invariant block}%
\end{equation}
are uniquely determined by $A$. The parameters $\sigma_{1},\ldots
,\sigma_{m_{2}}$ are the positive singular values of $V_{1}^{\ast}%
A\overline{V_{2}} $, so they are also uniquely determined by $A$.\medskip%
\newline
(b) $A$ is unitarily *congruent to a reduced form (\ref{Reduced Form}) in
which $m_{2}=m_{1}-\dim\left(  N(A)\cap N(A^{\ast})\right)  $; if $m_{2}>0$
then $D\in M_{m_{2}}$, $\Sigma=\operatorname{diag}(\sigma_{1},\ldots
,\sigma_{m_{2}})$, and all $\sigma_{i}>0$; if $m_{1}+m_{2}<n$, then
$A^{\prime}$ is square and $[A^{\prime}~B]$ has linearly independent rows; the
integers $m_{1}$, $m_{2}$ and the unitary *congruence class of the block
(\ref{invariant block}) are uniquely determined by $A$. The parameters
$\sigma_{1},\ldots,\sigma_{m_{2}}$ are the positive singular values of
$V_{1}^{\ast}AV_{2}$, so they are also uniquely determined by $A$.
\end{theorem}

\begin{proof}
We have%
\[
V^{\ast}A=\left[
\begin{array}
[c]{c}%
V_{1}^{\ast}A\\
V_{2}^{\ast}A
\end{array}
\right]  =\left[
\begin{array}
[c]{c}%
V_{1}^{\ast}A\\
0
\end{array}
\right]  \text{.}%
\]
The next step depends on whether we want to perform a unitary congruence or a
unitary *congruence. \medskip%
\newline
(a) Let $N=V_{1}^{\ast}A\overline{V_{2}}$ and form the unitary congruence%
\[
V^{\ast}A\left(  V^{\ast}\right)  ^{T}=\left[
\begin{array}
[c]{c}%
V_{1}^{\ast}A\\
0
\end{array}
\right]  \left[
\begin{array}
[c]{cc}%
\overline{V_{1}} & \overline{V_{2}}%
\end{array}
\right]  =\left[
\begin{array}
[c]{cc}%
V_{1}^{\ast}A\overline{V_{1}} & V_{1}^{\ast}A\overline{V_{2}}\\
0 & 0_{m_{1}}%
\end{array}
\right]  =\left[
\begin{array}
[c]{cc}%
M & N\\
0 & 0
\end{array}
\right]  \text{.}%
\]
Now let $m_{2}=\operatorname{rank}N$. If $m_{2}=0$ then $N=0$ and the form
(\ref{Reduced Form}) has been achieved with $A^{\prime}=M$. If $m_{2}>0$, use
the singular value decomposition to write $N=X\Sigma_{2}Y^{\ast}$, in which
$X$ and $Y$ are unitary,
\begin{equation}
\Sigma_{2}=\left[
\begin{array}
[c]{cc}%
0 & 0\\
\Sigma & 0
\end{array}
\right]  \text{ and }\Sigma=\operatorname{diag}(\sigma_{1},\ldots
,\sigma_{m_{2}})\text{,}\label{Sigma2}%
\end{equation}
and the diagonal entries $\sigma_{i}$ are the positive singular values of $N
$. Let $Z=X^{\ast}\oplus Y^{T}$ and form the unitary congruence%
\[
Z\left(  V^{\ast}A\bar{V}\right)  Z^{T}=\left[
\begin{array}
[c]{cc}%
X^{\ast}M\bar{X} & X^{\ast}NY\\
0 & 0_{m_{1}}%
\end{array}
\right]  =\left[
\begin{array}
[c]{cc|c}%
A^{\prime} & B & 0\\
C & D & \left[  \Sigma~0\right] \\\hline
\multicolumn{2}{c|}{0} & 0_{m_{1}}%
\end{array}
\right]  \hspace{-0.13in}%
\begin{array}
[c]{l}%
\\
\}m_{2}\\
\}m_{1}%
\end{array}
\text{.}%
\]
The block $X^{\ast}M\bar{X}$ has been partitioned so that $D\in M_{m_{2}}$.
Finally, inspection of (\ref{Reduced Form}) shows that $\dim(N(A)\cap
N(A^{T}))=m_{1}-m_{2}$.

Suppose that $R,\underline{R},U\in M_{n}$, $U$ is unitary, and $R=U\underline
{R}U^{T}$, so $R$ and $\underline{R}$ have the same parameters $m_{1}$ and
$m_{2}$. Suppose%
\[
R=\left[
\begin{array}
[c]{cc|c}%
A^{\prime} & B & 0\\
C & D & \left[  \Sigma~0\right] \\\hline
\multicolumn{2}{c|}{0} & 0_{m_{1}}%
\end{array}
\right]  \hspace{-0.13in}%
\begin{array}
[c]{l}%
\\
\}m_{2}\\
\}m_{1}%
\end{array}
\text{ and }\underline{R}=\left[
\begin{array}
[c]{cc|c}%
\underline{A}^{\prime} & \underline{B} & 0\\
\underline{C} & \underline{D} & \left[  \underline{\Sigma}~0\right] \\\hline
\multicolumn{2}{c|}{0} & 0_{m_{1}}%
\end{array}
\right]  \hspace{-0.13in}%
\begin{array}
[c]{l}%
\\
\}m_{2}\\
\}m_{1}%
\end{array}
\text{,}%
\]
partition $U=[U_{ij}]_{i,j=1}^{2}$ so that $U_{11}\in M_{n-m_{1}}$ and
$U_{22}\in M_{m_{1}}$, and partition%
\[
\underline{R}=\left[
\begin{array}
[c]{c}%
Z\\
0
\end{array}
\right]  \hspace{-0.13in}%
\begin{array}
[c]{l}%
\\
\}m_{1}%
\end{array}
\]
in which
\[
Z=\left[
\begin{array}
[c]{ccc}%
\underline{A}^{\prime} & \underline{B} & 0\\
\underline{C} & \underline{D} & \left[  \underline{\Sigma}~0\right]
\end{array}
\right]
\]
has full row rank. Then%
\[
\left[
\begin{array}
[c]{c}%
\bigstar\\
0
\end{array}
\right]  =R\bar{U}=U\underline{R}=\left[
\begin{array}
[c]{c}%
\bigstar\\
U_{21}Z
\end{array}
\right]  \text{,}%
\]
so $U_{21}=0$. Lemma \ref{Zero Blocks}(a) ensures that $U_{12}=0$ as well, so
$U=U_{11}\oplus U_{22}$ and both direct summands are unitary. Then%
\begin{align*}
R  & =\left[
\begin{array}
[c]{cc}%
\left[
\begin{array}
[c]{cc}%
A^{\prime} & B\\
C & D
\end{array}
\right]  & \left[
\begin{array}
[c]{c}%
0\\
\left[  \Sigma~0\right]
\end{array}
\right] \\
0 & 0_{m_{1}}%
\end{array}
\right] \\
& =U\underline{R}U^{T}=\left[
\begin{array}
[c]{cc}%
U_{11}\left[
\begin{array}
[c]{cc}%
\underline{A}^{\prime} & \underline{B}\\
\underline{C} & \underline{D}%
\end{array}
\right]  U_{11}^{T} & U_{11}\left[
\begin{array}
[c]{c}%
0\\
\left[  \underline{\Sigma}~0\right]
\end{array}
\right]  U_{22}^{T}\\
0 & 0_{m_{1}}%
\end{array}
\right]
\end{align*}
and the uniqueness assertion follows. \medskip%
\newline
(b) Let $N=V_{1}^{\ast}AV_{2}$ and form the unitary $^{\ast}$congruence%
\[
V^{\ast}A\left(  V^{\ast}\right)  ^{\ast}=\left[
\begin{array}
[c]{c}%
V_{1}^{\ast}A\\
V_{2}^{\ast}A
\end{array}
\right]  \left[
\begin{array}
[c]{cc}%
V_{1} & V_{2}%
\end{array}
\right]  =\left[
\begin{array}
[c]{cc}%
V_{1}^{\ast}AV_{1} & V_{1}^{\ast}AV_{2}\\
0 & 0_{m_{1}}%
\end{array}
\right]  =\left[
\begin{array}
[c]{cc}%
M & N\\
0 & 0
\end{array}
\right]  \text{.}%
\]
Let $m_{2}=\operatorname{rank}N$. If $m_{2}=0$ then $N=0$ and the form
(\ref{Reduced Form}) has been achieved with $A^{\prime}=M$. If $m_{2}>0$, use
the singular value decomposition to write $N=X\Sigma_{2}Y^{\ast}$, in which
$X$ and $Y$ are unitary, $\Sigma_{2}$ has the form (\ref{Sigma2}), and the
diagonal entries $\sigma_{i}$ are the positive singular values of $N$. Let
$Z=X^{\ast}\oplus Y^{\ast}$ and form the unitary *congruence%
\[
Z\left(  V^{\ast}AV\right)  Z^{\ast}=\left[
\begin{array}
[c]{cc}%
X^{\ast}MX & X^{\ast}NY\\
0 & 0_{m_{1}}%
\end{array}
\right]  =\left[
\begin{array}
[c]{cc|c}%
A^{\prime} & B & 0\\
C & D & \left[  \Sigma~0\right] \\\hline
\multicolumn{2}{c|}{0} & 0_{m_{1}}%
\end{array}
\right]  \hspace{-0.13in}%
\begin{array}
[c]{l}%
\\
\}m_{2}\\
\}m_{1}%
\end{array}
\text{.}%
\]
The block $X^{\ast}MX$ has been partitioned so that $D\in M_{m_{2}}$. The
uniqueness assertion follows from an argument parallel to the one employed in
(a).\hfill
\end{proof}

\medskip We are concerned here with only the simplest cases of unitary
congruence and *congruence, and the preceding theorem suffices for our
purpose; a general sparse form that can be achieved via unitary congruence and
*congruence is given in \cite[Theorem 6(d)]{HSregularization}.

\begin{corollary}
\label{Simplest Regularizations}Let $A\in M_{n}$ be singular and nonzero. Let
$m_{1}$ be the nullity of $A$, let the columns of $V_{1}$ be an orthonormal
basis for the range of $A$, let the columns of $V_{2}$ be an orthonormal basis
for the null space of $A^{\ast}$, and form the unitary matrix $V=[V_{1}%
~V_{2}]$. \medskip%
\newline
(a) Suppose $A$ is congruence normal. Then it is unitarily congruent to a
direct sum of the form%
\begin{equation}
A^{\prime}\oplus%
{\displaystyle\bigoplus\limits_{i=1}^{m_{2}}}
\sigma_{i}\left[
\begin{array}
[c]{cc}%
0 & 1\\
0 & 0
\end{array}
\right]  \oplus0_{m_{1}-m_{2}}\text{,\quad}\sigma_{i}>0\text{,}%
\label{regularization2}%
\end{equation}
in which $A^{\prime}$ is either absent or it is nonsingular and congruence
normal; $m_{2}=\operatorname{rank}A-\operatorname{rank}\bar{A}%
A=\operatorname{rank}V_{1}^{\ast}A\overline{V_{2}}$; and the parameters
$\sigma_{1},\ldots,\sigma_{m_{2}}$ are the positive singular values of
$V_{1}^{\ast}A\overline{V_{2}}$. The unitary congruence class of $A^{\prime}$,
$m_{2}$, $\sigma_{1}$,$\ldots$, $\sigma_{m_{2}}$, and $m_{1}$ are uniquely
determined by $A$. \medskip%
\newline
(b) Suppose $A$ is squared normal. Then it is unitarily *congruent to a direct
sum of the form (\ref{regularization2}), in which $A^{\prime}$ is either
absent or it is nonsingular and squared normal; $m_{2}=\operatorname{rank}%
A-\operatorname{rank}A^{2}=\operatorname{rank}V_{1}^{\ast}AV_{2}$; and the
parameters $\sigma_{1},\ldots,\sigma_{m_{2}}$ are the positive singular values
of $V_{1}^{\ast}AV_{2}$. The unitary *congruence class of $A^{\prime}$,
$m_{2}$, $\sigma_{1}$,$\ldots$, $\sigma_{m_{2}}$, and $m_{1}$ are uniquely
determined by $A$.
\end{corollary}

\begin{proof}
Combine Lemma \ref{Zero Blocks}(b) with Theorem \ref{UnitaryRegularization}%
.\hfill
\end{proof}

\medskip

The matrix $A^{\prime}$ in (\ref{regularization2}) is the \emph{regular part}
of $A$ under unitary congruence (respectively, unitary *congruence); the
direct sum of the singular summands in (\ref{regularization2}) is the
\emph{singular part} of $A$ under unitary congruence (respectively, unitary *congruence).

\section{Canonical forms}

We have now completed all the steps required to establish canonical forms for
a conjugate normal matrix $A$ under unitary congruence, and a squared normal
matrix $A$ under unitary *congruence: First apply the unitary regularization
described in Corollary \ref{Simplest Regularizations} to obtain the regular
and singular parts of $A$, then use Theorems \ref{nonsingularTheorem} and
\ref{*nonsingularTheorem} to identify the canonical form of the regular part.

\begin{theorem}
\label{General Congruence Normal}Let $A\in M_{n}$. If $\bar{A}A$ is normal,
then $A$ is unitarily congruent to a direct sum of blocks, each of which has
the form%
\begin{equation}
\left[  \sigma\right]  \text{ or }\tau\left[
\begin{array}
[c]{cc}%
0 & 1\\
\mu & 0
\end{array}
\right]  \text{, }\sigma,\tau\in\mathbb{R}\text{, }\sigma\geq0\text{, }%
\tau>0\text{, }\mu\in\mathbb{C}\text{, and }\mu\neq1\text{.}\label{gcon0}%
\end{equation}
This direct sum is uniquely determined by $A$ up to permutation of its blocks
and replacement of any parameter $\mu$ by $\mu^{-1}$. Conversely, if $A$ is
unitarily congruent to a direct sum of blocks of the form (\ref{gcon0}), then
$\bar{A}A$ is normal.
\end{theorem}

\begin{proof}
The unitary congruence regularization (\ref{regularization2}) reveals two
types of singular blocks
\begin{equation}
\lbrack0]\text{ and }\gamma\left[
\begin{array}
[c]{cc}%
0 & 1\\
0 & 0
\end{array}
\right]  \text{, }\gamma>0\text{,}\label{gcon1}%
\end{equation}
while Theorem \ref{nonsingularTheorem} reveals two types of nonsingular blocks%
\begin{equation}
\left[  \sigma\right]  \text{ and }\tau\left[
\begin{array}
[c]{cc}%
0 & 1\\
\mu & 0
\end{array}
\right]  \text{, }\sigma>0\text{, }\tau>0\text{, and }0\neq\mu\neq
1\text{.}\label{gcon2}%
\end{equation}

\hfill
\end{proof}

\begin{theorem}
\label{squared normal canonical 1}Let $A\in M_{n}$. If $A^{2}$ is normal, then
$A$ is unitarily *congruent to a direct sum of blocks, each of which is%
\begin{equation}
\lbrack\lambda]\text{ or }\tau\left[
\begin{array}
[c]{cc}%
0 & 1\\
\mu & 0
\end{array}
\right]  \text{, }\tau\in\mathbb{R}\text{, }\lambda,\mu\in\mathbb{C}\text{,
}\tau>0\text{, and }\left\vert \mu\right\vert <1\text{.}\label{g*0}%
\end{equation}
This direct sum is uniquely determined by $A$, up to permutation of its
blocks. Conversely, if $A$ is unitarily *congruent to a direct sum of blocks
of the form (\ref{g*0}), then $A^{2}$ is normal.
\end{theorem}

\begin{proof}
The unitary *congruence regularization (\ref{Reduced Form}) reveals the
singular blocks and Theorem \ref{*nonsingularTheorem} reveals the nonsingular
blocks.\hfill
\end{proof}

\medskip

For some applications, it can be convenient to know that the $2$-by-$2$ blocks
in (\ref{g*0}) may be replaced by canonical upper triangular blocks. The set
\begin{equation}
\mathcal{D}_{+}:=\{z\in\mathbb{C}:\operatorname{Re}z>0\}\cup\{it:t\in
\mathbb{R}\text{ and }t\geq0\}\label{D+}%
\end{equation}
has the useful property that every complex number has a unique square root in
$\mathcal{D}_{+}$. We use the following criterion of Pearcy:

\begin{lemma}
[\cite{PearcyPaper}]\label{Pearcy}Let $X,Y\in M_{2}$. Then $X$ and $Y$ are
unitarily *congruent if and only if $\operatorname{tr}X=\operatorname{tr}Y$,
$\operatorname{tr}X^{2}=\operatorname{tr}Y^{2}$, and $\operatorname{tr}%
X^{\ast}X=\operatorname{tr}Y^{\ast}Y$.
\end{lemma}

\begin{theorem}
\label{squared normal canonical 2}Let $A\in M_{n}$. If $A^{2}$ is normal, then
$A$ is unitarily *congruent to a direct sum of blocks, each of which is%
\begin{equation}
\lbrack\lambda]\text{ or }\left[
\begin{array}
[c]{cc}%
\nu & r\\
0 & -\nu
\end{array}
\right]  \text{, }\lambda,\nu\in\mathbb{C}\text{, }r\in\mathbb{R}\text{,
}r>0\text{, and }\nu\in\mathcal{D}_{+}\text{.}\label{t*0}%
\end{equation}
This direct sum is uniquely determined by $A$ up to permutation of its blocks.
Conversely, if $A$ is unitarily *congruent to a direct sum of blocks of the
form (\ref{t*0}), then $A^{2}$ is normal.
\end{theorem}

\begin{proof}
It suffices to show that if $\tau>0$ and $\left\vert \mu\right\vert <1$, and
if we define
\begin{equation}
\nu:=\tau\sqrt{\mu}\in\mathcal{D}_{+}\label{t*1}%
\end{equation}
and%
\begin{equation}
r:=\tau\left(  1-\left\vert \mu\right\vert \right)  \text{,}\label{t*2}%
\end{equation}
then%
\[
C_{1}:=\left[
\begin{array}
[c]{cc}%
\nu & r\\
0 & -\nu
\end{array}
\right]  \text{ and }C_{2}:=\tau\left[
\begin{array}
[c]{cc}%
0 & 1\\
\mu & 0
\end{array}
\right]
\]
are unitarily *congruent. One checks that%
\[
\operatorname{tr}C_{1}=0=\operatorname{tr}C_{2}\text{,}%
\]%
\[
\operatorname{tr}C_{1}^{2}=2\nu^{2}=2\tau^{2}\mu=\operatorname{tr}C_{2}%
^{2}\text{,}%
\]
and%
\[
\operatorname{tr}C_{1}^{\ast}C_{1}=2\left\vert \nu\right\vert ^{2}+r^{2}%
=2\tau^{2}\left\vert \mu\right\vert +\tau^{2}\left(  1-\left\vert
\mu\right\vert \right)  ^{2}=\tau^{2}+\tau^{2}\mu^{2}=\operatorname{tr}%
C_{2}^{\ast}C_{2}\text{,}%
\]
so our assertion follows from Lemma \ref{Pearcy}.\hfill
\end{proof}

\section{Beyond normality}

We conclude with several results involving unitary congruence and unitary *congruence.

\subsection{Criteria for unitary congruence and *congruence}

To show that two matrices are unitarily congruent (or unitarily *congruent),
in certain cases it suffices to show only that they are congruent (or *congruent).

\begin{theorem}
\label{Criterion}Let $A,B,S\in M_{n}$ be nonsingular and let $S=WQ$ be a right
polar decomposition. Then \medskip%
\newline
(a) $A$ and $B$ are unitarily congruent if and only if the pairs $(A,A^{-\ast
})$ and $(B,B^{-\ast})$ are congruent. In fact, if $(A,A^{-\ast}%
)=S(B,B^{-\ast})S^{T}$, then $(A,A^{-\ast})=W(B,B^{-\ast})W^{T}$. \medskip%
\newline
(b) $A$ and $B$ are unitarily *congruent if and only if the pairs
$(A,A^{-\ast})$ and $(B,B^{-\ast})$ are *congruent. In fact, if $(A,A^{-\ast
})=S(B,B^{-\ast})S^{\ast}$, then $(A,A^{-\ast})=W(B,B^{-\ast})W^{\ast}$.
\end{theorem}

\begin{proof}
Let $\lambda_{1}>\cdots>\lambda_{d}>0$ be the distinct eigenvalues of
$S^{\ast}S$ and let $p(t)$ be any polynomial such that $p(\lambda
_{i})=+\lambda_{i}^{1/2}$ and $p(\lambda_{i}^{-1})=+\lambda_{i}^{-1/2}$ for
each $i=1,\ldots,d$. Then $Q=p(S^{\ast}S)$ is Hermitian and positive definite,
$Q^{2}=S^{\ast}S$, and $Q^{-1}=p\left(  (S^{\ast}S)^{-1}\right)  $.
\newline
(a) If there is a unitary $U$ such that $A=UBU^{T}$, then
\[
(A,A^{-\ast})=(UBU^{T},UB^{-\ast}U^{T})=U(B,B^{-\ast})U^{T}\text{.}%
\]
Conversely, if $(A,A^{-\ast})=S(B,B^{-\ast})S^{T}$, then%
\[
SBS^{T}=A=\left(  A^{-\ast}\right)  ^{-\ast}=\left(  SB^{-\ast}S^{T}\right)
^{-\ast}=S^{-\ast}B\bar{S}^{-1}\text{,}%
\]
so%
\[
\left(  S^{\ast}S\right)  B=B\left(  S^{\ast}S\right)  ^{-T}\text{.}%
\]
It follows that%
\[
g(S^{\ast}S)B=Bg(\left(  S^{\ast}S\right)  ^{-1})^{T}%
\]
for any polynomial $g(t)$. Choosing $g(t)=p(t)$, we have%
\[
QB=p(S^{\ast}S)B=Bp(\left(  S^{\ast}S\right)  ^{-1})^{T}=BQ^{-T}\text{,}%
\]
so $QBQ^{T}=B$ and%
\[
A=SBS^{T}=WQBQ^{T}W^{T}=WBW^{T}\text{.}%
\]
(b) If there is a unitary $U$ such that $A=UBU^{\ast}$, then
\[
(A,A^{-\ast})=(UBU^{T},UB^{-\ast}U^{\ast})=U(B,B^{-\ast})U^{\ast}\text{.}%
\]
Conversely, if $(A,A^{-\ast})=S(B,B^{-\ast})S^{\ast}$, then%
\[
SBS^{\ast}=A=\left(  A^{-\ast}\right)  ^{-\ast}=\left(  SB^{-\ast}S^{\ast
}\right)  ^{-\ast}=S^{-\ast}BS^{-1}\text{,}%
\]
so%
\[
\left(  S^{\ast}S\right)  B=B\left(  S^{\ast}S\right)  ^{-1}%
\]
and hence%
\[
g(S^{\ast}S)B=Bg(\left(  S^{\ast}S\right)  ^{-1})
\]
for any polynomial $g(t)$. Choosing $g(t)=p(t)$, we have%
\[
QB=p(S^{\ast}S)B=Bp(\left(  S^{\ast}S\right)  ^{-1})=BQ^{-1}\text{,}%
\]
so $QBQ=QBQ^{\ast}=B$ and%
\[
A=SBS^{\ast}=WQBQ^{\ast}W^{\ast}=WQBQW^{\ast}=WBW^{\ast}\text{.}%
\]

\hfill
\end{proof}

\begin{corollary}
\label{CriterionCorollary}Let $A,B\in M_{n}$ be given. \medskip%
\newline
(a) If $A$ and $B$ are either both unitary or both coninvolutory, then $A$ and
$B$ are unitarily congruent if and only if they are congruent. \medskip%
\newline
(b) If $A$ and $B$ are either both unitary or both involutory, then $A$ and
$B$ are unitarily *congruent if and only if they are *congruent.
\end{corollary}

\begin{proof}
The key observation is that $A^{-\ast}=A$ if $A$ is unitary, $A^{-\ast}=A^{T}$
if $A$ is coninvolutory, and $A^{-\ast}=A^{\ast}$ if $A$ is involutory.
\medskip%
\newline
(a) Suppose $A=SBS^{T}$. If $A$ and $B$ are unitary, then%
\[
(A,A^{-\ast})=(A,A)=(SBS^{T},SBS^{T})=S(B,B)S^{T}=S(B,B^{-\ast})S^{T}\text{,}%
\]
so Theorem \ref{Criterion}(a) ensures that $A$ is unitarily congruent to $B$.
If $A$ and $B$ are coninvolutory, then%
\[
(A,A^{-\ast})=(A,A^{T})=(SBS^{T},SB^{T}S^{T})=S(B,B^{T})S^{T}=S(B,B^{-\ast
})S^{T}\text{,}%
\]
so $A$ is again unitarily congruent to $B$. \medskip%
\newline
(b) Suppose $A=SBS^{\ast}$. If $A$ and $B$ are unitary, then
\[
(A,A^{-\ast})=(A,A)=(SBS^{\ast},SBS^{\ast})=S(B,B)S^{\ast}=S(B,B^{-\ast
})S^{\ast}\text{,}%
\]
so $A$ is unitarily *congruent to $B$. If $A$ and $B$ are involutory, then%
\[
(A,A^{-\ast})=(A,A^{\ast})=(SBS^{\ast},SB^{\ast}S^{\ast})=S(B,B^{\ast}%
)S^{\ast}=S(B,B^{-\ast})S^{\ast}\text{,}%
\]
so $A$ is unitarily congruent to $B$.\hfill
\end{proof}

\subsection{Hermitian cosquares}

\begin{theorem}
\label{HermitianCosquaresTheorem}Suppose that $A\in M_{n}$ is nonsingular. The
following are equivalent:\medskip%
\newline
(a) $A^{-T}A$ is Hermitian. \medskip%
\newline
(b) $\bar{A}A$ is Hermitian. \medskip%
\newline
(c) $A$ is unitarily congruent to a direct sum of real blocks, each of which
is%
\begin{equation}
\left[  \sigma\right]  \text{ or }\tau\left[
\begin{array}
[c]{cc}%
0 & 1\\
\mu & 0
\end{array}
\right]  \text{,\quad}\sigma>0\text{, }\tau>0\text{, }\mu\in\mathbb{R}\text{,
}0\neq\mu\neq1\text{.}\label{HermitianCongruence}%
\end{equation}
This direct sum is uniquely determined by $A$, up to permutation of its blocks
and replacement of any $\mu$ by $\mu^{-1}$. Conversely, if $A$ is unitarily
congruent to a direct sum of blocks of the form (\ref{HermitianCongruence}),
then $\bar{A}A$ is Hermitian.
\end{theorem}

\begin{proof}
$A^{-T}A$ is Hermitian if and only if%
\[
A^{-T}A=\left(  A^{-T}A\right)  ^{\ast}=A^{\ast}\bar{A}^{-1}%
\]
if and only if%
\[
A\bar{A}=A^{T}A^{\ast}=\left(  A\bar{A}\right)  ^{\ast}%
\]
if and only if $\bar{A}A=\overline{A\bar{A}}$ is Hermitian. The canonical
blocks (\ref{HermitianCongruence}) are the same as those in (\ref{cnc0}), with
the restriction that $\mu$ must be real.\hfill
\end{proof}

\medskip

Any coninvolution $A$ satisfies the hypotheses of the preceding theorem:
$\bar{A}A=I$.

\begin{corollary}
Suppose that $A\in M_{n}$ and $\bar{A}A=I$. Then $A$ is unitarily congruent to%
\begin{equation}
I_{n-2q}\oplus%
{\displaystyle\bigoplus\limits_{j=1}^{q}}
\left[
\begin{array}
[c]{cc}%
0 & \sigma_{j}^{-1}\\
\sigma_{j} & 0
\end{array}
\right]  \text{,\quad}\sigma_{j}>1\text{,}\label{ConinvolutionCongruence}%
\end{equation}
in which $\sigma_{1},\sigma_{1}^{-1},\ldots,\sigma_{q},\sigma_{q}^{-1}$ are
the singular values of $A$ that are different from $1$ and each $\sigma_{j}%
>1$. Conversely, if $A$ is unitarily congruent to a direct sum of the form
(\ref{ConinvolutionCongruence}), then $A$ is coninvolutory. Two coninvolutions
of the same size are unitarily congruent if and only if they have the same
singular values.
\end{corollary}

\begin{proof}
$\bar{A}A$ is Hermitian, so $A$ is unitarily congruent to a direct sum of
blocks of the two types (\ref{HermitianCongruence}). But $\bar{A}A=I$, so
$\sigma=1$ and $\tau^{2}\mu=1$. Then $\tau=\left(  \tau\mu\right)  ^{-1}$, so%
\[
\tau\left[
\begin{array}
[c]{cc}%
0 & 1\\
\mu & 0
\end{array}
\right]  =\left[
\begin{array}
[c]{cc}%
0 & \left(  \tau\mu\right)  ^{-1}\\
\tau\mu & 0
\end{array}
\right]  \text{,}%
\]
which has singular values $\tau\mu$ and $\left(  \tau\mu\right)  ^{-1}$.\hfill
\end{proof}

\medskip

The general case is obtained by specializing Theorem
\ref{General Congruence Normal}.

\begin{theorem}
Let $A\in M_{n}$ and suppose that $\bar{A}A$ is Hermitian. Then $A$ is
unitarily congruent to a direct sum of blocks, each of which is%
\begin{equation}
\left[  \sigma\right]  \text{ or }\tau\left[
\begin{array}
[c]{cc}%
0 & 1\\
\mu & 0
\end{array}
\right]  \text{, }\sigma,\tau,\mu\in\mathbb{R}\text{, }\sigma\geq0\text{, }%
\mu\neq1\text{.}\label{HermitianUnitaryCongruenceBlocks}%
\end{equation}
This direct sum is uniquely determined by $A$ up to permutation of its blocks
and replacement of any (real) parameter $\mu$ by $\mu^{-1}$. Conversely, if
$A$ is unitarily congruent to a direct sum of blocks of the form
(\ref{HermitianUnitaryCongruenceBlocks}), then $\bar{A}A$ is Hermitian.
\end{theorem}

\subsection{Unitary cosquares}

\begin{theorem}
\label{NonsingularConjugateNormal}Suppose that $A\in M_{n}$ is nonsingular.
The following are equivalent:\medskip%
\newline
(a) $A^{-T}A$ is unitary.\medskip%
\newline
(b) $A$ is conjugate normal.\medskip%
\newline
(c) $A$ is unitarily congruent to a direct sum of blocks, each of which is%
\begin{equation}
\left[  \sigma\right]  \text{ or }\tau\left[
\begin{array}
[c]{cc}%
0 & 1\\
e^{i\theta} & 0
\end{array}
\right]  \text{, }\sigma,\tau,\theta\in\mathbb{R}\text{, }\sigma>0\text{,
}\tau>0\text{, }0<\theta\leq\pi\text{.}\label{ConjugateNormalCongruence}%
\end{equation}
This direct sum is uniquely determined by the eigenvalues of $\bar{A}A$, up to
permutation of its summands. If $A$ is unitarily congruent to a direct sum of
blocks of the form (\ref{ConjugateNormalCongruence}), then $A$ is conjugate normal.
\end{theorem}

\begin{proof}
$A^{-T}A$ is unitary if and only if%
\[
A^{-1}A^{T}=\left(  A^{-T}A\right)  ^{-1}=\left(  A^{-T}A\right)  ^{\ast
}=A^{\ast}\bar{A}^{-1}%
\]
if and only if%
\[
AA^{\ast}=A^{T}\bar{A}=\overline{A^{\ast}A}\text{.}%
\]
The canonical blocks (\ref{ConjugateNormalCongruence}) follow from
(\ref{cnc0}) by specialization.

The eigenvalues of $\bar{A}A$ are of two types: positive eigenvalues that
correspond to squares of blocks of the first type in
(\ref{ConjugateNormalCongruence}), and conjugate pairs of non-positive (but
possibly negative) eigenvalues $\{\tau^{2}e^{i\theta},\tau^{2}e^{-i\theta}\}$
that correspond to blocks of the second type with $0<\theta\leq\pi$. Thus, the
parameters $\sigma$, $\tau$, and $e^{i\theta}$ of the blocks in
(\ref{ConjugateNormalCongruence}) can be inferred from the eigenvalues of
$\bar{A}A$.\hfill
\end{proof}

\medskip

The unitary congruence canonical blocks (\ref{ConjugateNormalCongruence}) for
a conjugate normal matrix are a subset of the canonical blocks (\ref{gcon0})
for a congruence normal matrix; the 2-by-2 singular blocks are omitted, and
the 2-by-2 nonsingular blocks are required to be positive scalar multiples of
a unitary block. This observation shows that every conjugate normal matrix is
congruence normal. Moreover, examination of the canonical blocks of a
conjugate normal matrix shows that it is unitarily congruent to a direct sum
of a positive diagonal matrix, positive scalar multiples of unitary matrices
(which the following corollary shows may be taken to be real), and a zero
matrix. Thus, a conjugate normal matrix is unitarily congruent to a real
normal matrix.

If $A$ itself is unitary, then its cosquare $A^{-T}A=\bar{A}A$ is certainly
unitary, so the unitary congruence canonical form of a unitary matrix follows
from the preceding theorem. Of course, the eigenvalues of the cosquare of a
unitary matrix are all unimodular and are constrained by the conditions in
(\ref{JCFdiagonalizableCosquare}): any eigenvalue $\mu\neq1$ (even $\mu=-1$)
occurs in a conjugate pair $\{\mu,\bar{\mu}\}$.

\begin{corollary}
\label{UnitarySpecialCase}Suppose that $U\in M_{n}$ is unitary. Then $U$ is
unitarily congruent to%
\begin{equation}
I_{n-2q}\oplus%
{\displaystyle\bigoplus\limits_{j=1}^{q}}
\left[
\begin{array}
[c]{cc}%
0 & 1\\
\mu_{j} & 0
\end{array}
\right]  \text{,\quad\ }\mu\in\mathbb{C}\text{, }\left\vert \mu_{j}\right\vert
=1\text{, }\mu_{j}\neq1\text{,}\label{UnitarySpecialCaseBlocks}%
\end{equation}
in which $\mu_{1},\bar{\mu}_{1},\ldots,\mu_{q},\bar{\mu}_{q}$ are the
eigenvalues of $\bar{U}U$ that are different from $1$. If $\mu_{j}%
=e^{i\theta_{j}}$, then each unitary 2-by-2 block $H_{2}(e^{i\theta_{j}})$ in
(\ref{UnitarySpecialCaseBlocks}) can be replaced by a real orthogonal block
\begin{equation}
Q_{2}(\theta)=\left[
\begin{array}
[c]{cc}%
\alpha & \beta\\
-\beta & \alpha
\end{array}
\right] \label{RealBlock}%
\end{equation}
in which $\alpha=\cos(\theta/2)$ and $\beta=\sin(\theta/2)$, or by a Hermitian
unitary block%
\begin{equation}
\mathcal{H}_{2}(\theta)=\left[
\begin{array}
[c]{cc}%
0 & e^{-i\theta/2}\\
e^{i\theta/2} & 0
\end{array}
\right] \label{HermitianBlock}%
\end{equation}
Thus, $U$ is unitarily congruent to a real orthogonal matrix as well as to a
Hermitian unitary matrix.
\end{corollary}

\begin{proof}
$U^{-T}U=\bar{U}U$ so the parameters $\mu$ in (\ref{ConjugateNormalCongruence}%
) correspond to the pairs of unimodular conjugate eigenvalues of $\bar{U}U$.
Since each block in (\ref{ConjugateNormalCongruence}) must be unitary, the
parameters $\sigma$ and $\tau$ must be $+1$. One checks that the cosquares of
$H_{2}(e^{i\theta_{j}})$ and $Q_{2}(\theta)$ (\ref{RealBlock}) (both unitary)
have the same eigenvalues (namely, $e^{\pm i\theta_{j}})$, so they are
similar. Theorem \ref{CongruenceCanonicalForms}(a) ensures that $H_{2}%
(e^{i\theta_{j}})$ and $Q_{2}(\theta)$ are congruent, and Corollary
\ref{CriterionCorollary}(a) ensures that they are actually unitarily
congruent. The unitary congruence%
\[
\left[
\begin{array}
[c]{cc}%
e^{-i\theta/4} & 0\\
0 & e^{-i\theta/4}%
\end{array}
\right]  \left[
\begin{array}
[c]{cc}%
0 & 1\\
e^{i\theta} & 0
\end{array}
\right]  \left[
\begin{array}
[c]{cc}%
e^{-i\theta/4} & 0\\
0 & e^{-i\theta/4}%
\end{array}
\right]  =\left[
\begin{array}
[c]{cc}%
0 & e^{-i\theta/2}\\
e^{i\theta/2} & 0
\end{array}
\right]
\]
shows that the 2-by-2 blocks in (\ref{UnitarySpecialCaseBlocks}) may be
replaced by Hermitian blocks of the form $\mathcal{H}_{2}(\theta)$.\hfill
\end{proof}

\medskip

In order to state the general case of Theorem \ref{NonsingularConjugateNormal}%
, we need to know what the singular part of a conjugate normal matrix is,
after regularization by unitary congruence.

\begin{lemma}
Let $A\in M_{n}$ be singular and conjugate normal; let $m_{1}$ be the nullity
of $A$. Then \medskip%
\newline
(a) The angle between $Ax$ and $Ay$ is the same as the angle between $A^{T}x$
and $A^{T}y$ for all $x,y\in\mathbb{C}^{n}$. \medskip%
\newline
(b) $\left\Vert Ax\right\Vert =\left\Vert A^{T}x\right\Vert $ for all
$x\in\mathbb{C}^{n}$; \medskip%
\newline
(c) $N(A)=N(A^{T})$; and \medskip%
\newline
(d) $A$ is unitarily congruent to $A^{\prime}\oplus0_{m_{1}}$ in which
$A^{\prime}$ is nonsingular and conjugate normal.
\end{lemma}

\begin{proof}
Compute
\[
(Ax)^{\ast}(Ay)=x^{\ast}A^{\ast}Ay=x^{\ast}\overline{AA^{\ast}}y=(A^{T}%
x)^{\ast}(A^{T}y)\text{;}%
\]
when $x=y$ we have $\left\Vert Ax\right\Vert ^{2}=\left\Vert A^{T}x\right\Vert
^{2}$. In particular, $Ax=0$ if and only if $A^{T}x=0$.

In the reduced form (\ref{Reduced Form}) of $A$ we have
\[
m_{2}=\dim N(A)-\dim(N(A)\cap N(A^{T}))=\dim N(A)-\dim N(A)=0\text{.}%
\]
Thus, $A$ is unitarily congruent to $A^{\prime}\oplus0_{m_{1}}$; $A^{\prime}$
is nonsingular and its unitary congruence class is uniquely determined; and%
\[
\left(  A^{\prime}\right)  ^{\ast}A^{\prime}\oplus0_{m_{1}}=A^{\ast
}A=\overline{AA^{\ast}}=\overline{A^{\prime}\left(  A^{\prime}\right)  ^{\ast
}}\oplus0_{m_{1}}\text{,}%
\]
so $A^{\prime}$ is conjugate normal.\hfill
\end{proof}

\medskip

\begin{corollary}
\label{ConjugateNormalUnitaryCongruenceCorollary}Let $A\in M_{n}$ and suppose
that $A$ is conjugate normal. Then $A$ is unitarily congruent to a direct sum
of blocks, each of which is%
\begin{equation}
\left[  \sigma\right]  \text{ or }\tau\left[
\begin{array}
[c]{cc}%
0 & 1\\
e^{i\theta} & 0
\end{array}
\right]  \text{, }\sigma,\tau,\theta\in\mathbb{R}\text{, }\sigma\geq0\text{,
}\tau>0\text{, }0<\theta\leq\pi\text{.}%
\label{ConjugateNormalUnitaryCongruenceBlocks}%
\end{equation}
This direct sum is uniquely determined by the eigenvalues of $\bar{A}A$, up to
permutation of its blocks: there is one block $\sqrt{\rho}H_{2}(e^{i\theta})$
(with $\sqrt{\rho}>0$) corresponding to each conjugate eigenvalue pair $\{\rho
e^{i\theta},\rho e^{-i\theta}\}$ of $\bar{A}A$ with $\rho>0$ and $0<\theta
\leq\pi$; the number of blocks $[\sigma]$ with $\sigma>0$ equals the
multiplicity of $\sigma$ as a (positive) eigenvalue of $\bar{A}A$, so the
total number of blocks of this type equals the number of positive eigenvalues
of $\bar{A}A$; the number of blocks $[0]$ equals the nullity of $A$. \medskip%
\newline
If $B\in M_{n}$ is conjugate normal, then $A$ is unitarily congruent to $B$ if
and only if $\bar{A}A$ and $\bar{B}B$ have the same eigenvalues. \medskip%
\newline
Each unitary block $H_{2}(e^{i\theta})$ in
(\ref{ConjugateNormalUnitaryCongruenceBlocks}) can be replaced by a real
orthogonal block%
\[
\left[
\begin{array}
[c]{cc}%
\alpha & \beta\\
-\beta & \alpha
\end{array}
\right]  \text{, }\alpha=\cos(\theta/2)\text{ and }\beta=\sin(\theta
/2)\text{.}%
\]
If $A$ is unitarily congruent to a direct sum of blocks of the form
(\ref{ConjugateNormalUnitaryCongruenceBlocks}), then $A$ is conjugate
normal.\hfill
\end{corollary}

\subsection{Hermitian *cosquares}

\begin{theorem}
\label{Hermitian*cosquaresTheorem}Suppose that $A\in M_{n}$ is nonsingular.
The following are equivalent:\medskip%
\newline
(a) $A^{-\ast}A$ is Hermitian. \medskip%
\newline
(b) $A^{2}$ is Hermitian. \medskip%
\newline
(c) $A$ is unitarily *congruent to a direct sum (uniquely determined by $A$ up
to permutation of summands) of blocks, each of which is%
\begin{equation}
\left[  \lambda\right]  \text{, }\left[  i\nu\right]  \text{, or }\tau\left[
\begin{array}
[c]{cc}%
0 & 1\\
\mu & 0
\end{array}
\right]  \text{,\quad}\lambda,\nu,\mu\in\mathbb{R}\text{, }\lambda\neq0\neq
\nu\text{, }\tau>0\text{, }0<\left\vert \mu\right\vert <1\text{.}%
\label{Hermitian*Congruence}%
\end{equation}
If $\mu_{1},\mu_{1}^{-1},\ldots,\mu_{q},\mu_{q}^{-1}$ are the (real)
eigenvalues of $A^{-\ast}A$ that are not equal to $\pm1$ and satisfy
$0<\left\vert \mu_{j}\right\vert <1$ for $j=1,\ldots,q$, and if $A^{-\ast}A $
has $p$ eigenvalues equal to $+1$, then the unitary *congruence canonical form
of $A$ has $p$ blocks of the first type in (\ref{Hermitian*Congruence}),
$n-2q-p$ blocks of the second type, and $q$ blocks of the third type.
\end{theorem}

\begin{proof}
$A^{-\ast}A$ is Hermitian if and only if%
\[
A^{-\ast}A=\left(  A^{-\ast}A\right)  ^{\ast}=A^{\ast}A^{-1}%
\]
if and only if%
\[
A^{2}=\left(  A^{2}\right)  ^{\ast}\text{.}%
\]
The canonical blocks (\ref{Hermitian*Congruence}) follow from
(\ref{*nonsingular}) by specialization.\hfill
\end{proof}

\medskip

Any involutory matrix satisfies the hypotheses of the preceding theorem.

\begin{corollary}
Let $A\in M_{n}$, suppose that $A^{2}=I$, and suppose that $A$ has $p$
eigenvalues equal to $1$. The singular values of $A$ that are different from
$1$ occur in reciprocal pairs: $\sigma_{1},\sigma_{1}^{-1},\ldots,\sigma
_{q},\sigma_{q}^{-1}$ in which each $\sigma_{i}>1$. Then $A$ is unitarily
*congruent to%
\begin{equation}
I_{p-q}\oplus\left(  -I_{n-p-q}\right)
{\displaystyle\bigoplus\limits_{j=1}^{q}}
\left[
\begin{array}
[c]{cc}%
0 & \sigma_{j}^{-1}\\
\sigma_{j} & 0
\end{array}
\right]  \text{,\quad}\sigma_{j}>1\label{Involution*Congruence}%
\end{equation}
as well as to%
\begin{equation}
I_{p-q}\oplus\left(  -I_{n-p-q}\right)
{\displaystyle\bigoplus\limits_{i=1}^{q}}
\left[
\begin{array}
[c]{cc}%
1 & \sigma_{i}-\sigma_{i}^{-1}\\
0 & -1
\end{array}
\right]  \text{.}\label{Involution*Congruence2}%
\end{equation}
Conversely, if $A$ is unitarily *congruent to a direct sum of either form
(\ref{Involution*Congruence}) or (\ref{Involution*Congruence2}), then $A$ is
an involution and has $p-q$ eigenvalues equal to $1$. Two involutions of the
same size are unitarily *congruent if and only if they have the same singular
values and $+1$ is an eigenvalue with the same multiplicity for each of them.
\end{corollary}

\begin{proof}
Since $A=A^{-1}$, $A^{-\ast}A=A^{\ast}A$ and the eigenvalues of the *cosquare
are just the squares of the singular values of $A$; the eigenvalues of the
*cosquare $A^{\ast}A$ that are not equal to $1$ must occur in reciprocal
pairs. Let $\sigma_{1},\ldots,\sigma_{q}$ be the singular values of $A$ that
are greater than $1$. Each 2-by-2 block in (\ref{Hermitian*Congruence}) has
the form%
\[
\tau_{j}\left[
\begin{array}
[c]{cc}%
0 & 1\\
\sigma_{j}^{2} & 0
\end{array}
\right]  \text{,}%
\]
which has singular values $\tau_{j}\sigma_{j}^{2}$ and $\tau_{j}$; they are
reciprocal if and only if $\tau_{j}=\sigma_{j}^{-1}$. Each 2-by-2 block
contributes a pair of eigenvalues $\pm1$, which results in the asserted
summands $I_{p-q}\oplus\left(  -I_{n-p-q}\right)  $ since all of the
eigenvalues of $A$ are $\pm1$.

To confirm that $A$ is unitarily *congruent to the direct sum
(\ref{Involution*Congruence2}), it suffices to show that
\[
C_{1}=\left[
\begin{array}
[c]{cc}%
1 & \sigma-\sigma^{-1}\\
0 & -1
\end{array}
\right]  \quad\text{and\quad}C_{2}=\left[
\begin{array}
[c]{cc}%
1 & \sigma^{-1}\\
\sigma & 0
\end{array}
\right]
\]
are unitarily *congruent. Using Lemma \ref{Pearcy}, it suffices to observe
that%
\begin{align*}
\operatorname{tr}C_{1}  & =0=\operatorname{tr}C_{2}\\
\operatorname{tr}C_{1}^{2}  & =2=\operatorname{tr}C_{2}^{2}\text{, and}\\
\operatorname{tr}C_{1}^{\ast}C_{1}  & =\sigma^{2}+\sigma^{-2}%
=\operatorname{tr}C_{2}^{\ast}C_{2}%
\end{align*}

\hfill
\end{proof}

\medskip

The general case is obtained by specialization of Theorems
\ref{squared normal canonical 1} and \ref{squared normal canonical 2}.

\begin{theorem}
Let $A\in M_{n}$ and suppose $A^{2}$ is Hermitian. Then $A$ is unitarily
*congruent to a direct sum of blocks, each of which is%
\begin{equation}
\left[  \lambda\right]  \text{, }\left[  i\lambda\right]  \text{, or }%
\tau\left[
\begin{array}
[c]{cc}%
0 & 1\\
\mu & 0
\end{array}
\right]  \text{, }\lambda,\mu,\tau\in\mathbb{R}\text{, }\tau>0\text{, }%
-1<\mu<1\text{.}\label{sH1}%
\end{equation}
Alternatively, $A$ is unitarily *congruent to a direct sum of blocks, each of
which is%
\[
\left[  \lambda\right]  \text{, }\left[  i\lambda\right]  \text{, or }\left[
\begin{array}
[c]{cc}%
\tau\sqrt{\mu} & \tau\left(  1-\left\vert \mu\right\vert \right) \\
0 & -\tau\sqrt{\mu}%
\end{array}
\right]  \text{, }\sqrt{\mu}\in\mathcal{D}_{+}\text{,}%
\]
in which the parameters $\tau$ and $\mu$ satisfy the conditions in (\ref{sH1}).
\end{theorem}

\subsection{Unitary *cosquares}

\begin{theorem}
\label{Unitary*cosquaresTheorem}Suppose $A\in M_{n}$ is nonsingular. The
following are equivalent:\medskip%
\newline
(a) $A^{-\ast}A$ is unitary. \medskip%
\newline
(b) $A$ is normal. \medskip%
\newline
(c) $A$ is unitarily *congruent to a direct sum of blocks, each of which is%
\begin{equation}
\left[  \lambda\right]  \text{,\quad}\lambda\in\mathbb{C}\text{, }\lambda
\neq0\text{.}\label{Normal*Congruence}%
\end{equation}

\end{theorem}

\begin{proof}
$A^{-\ast}A$ is unitary if and only if%
\[
A^{-1}A^{\ast}=\left(  A^{-\ast}A\right)  ^{-1}=\left(  A^{-\ast}A\right)
^{\ast}=A^{\ast}A^{-1}%
\]
if and only if%
\[
AA^{\ast}=A^{\ast}A\text{.}%
\]
The canonical blocks (\ref{Normal*Congruence}) follow from (\ref{*nonsingular}%
) by specialization.\hfill
\end{proof}

\subsection{Projections and $\lambda$-projections}

The unitary *congruence regularization algorithm described in Theorem
\ref{UnitaryRegularization}(b) permits us to identify a unitary *congruence
canonical form for $\lambda$-projections, that is, matrices $A\in M_{n}$ such
that $A^{2}=\lambda A$. A $1$-projection is an ordinary projection, while a
nonzero $0$-projection is a nilpotent matrix with index $2$. A complex matrix
whose minimal polynomial is quadratic is a translation of a $\lambda$-projection.

\begin{theorem}
\label{LambdaProjection}Let $A\in M_{n}$ be singular and nonzero, let
$\lambda\in\mathbb{C}$ be given, and suppose that $A^{2}=\lambda A$. Let
$m_{1}$ be the nullity of $A$ and let $\tau_{1},\ldots,\tau_{m_{2}}$ be the
singular values of $A$ that are strictly greater than $|\lambda|$ ($m_{2}=0$
is possible). Then \medskip%
\newline
(a) $A$ is unitarily *congruent to
\begin{equation}
\lambda I_{n-m_{1}-m_{2}}\oplus%
{\displaystyle\bigoplus\limits_{i=1}^{m_{2}}}
\left[
\begin{array}
[c]{cc}%
\lambda & \sqrt{\tau_{i}^{2}-\left\vert \lambda\right\vert ^{2}}\\
0 & 0
\end{array}
\right]  \oplus0_{m_{1}-m_{2}}\text{.}\label{ProjectorBlocks}%
\end{equation}
This direct sum is uniquely determined by $\lambda$ and the singular values of
$A$, which are $\tau_{1},\ldots,\tau_{m_{2}}$, $\left\vert \lambda\right\vert
$ with multiplicity $n-m_{1}-m_{2}$, and $0$ with multiplicity $m_{1}$.
\medskip%
\newline
(b) For a given $\lambda$, two $\lambda$-projections of the same size are
unitarily *congruent if and only if they have the same rank and the same
singular values. \medskip%
\newline
(c) Suppose $0\neq A\neq\lambda I$ and let $\nu=\min\{m_{1},n-m_{1}\}$. Then
$\nu>0$ and the $\nu$ largest singular values of $A$ and $A-\lambda I$ are the
same. In particular, the spectral norms of $A$ and $A-\lambda I$ are equal.
\end{theorem}

\begin{proof}
(a) Let $F$ denote a reduced form (\ref{Reduced Form}) for $A$ under unitary
*congruence, which is also a $\lambda$-projection. Compute%
\[
F^{2}=\left[
\begin{array}
[c]{ccc}%
A^{\prime} & B & 0\\
C & D & \left[  \Sigma~0\right] \\
0 & 0 & 0_{m_{1}}%
\end{array}
\right]  ^{2}=\left[
\begin{array}
[c]{ccc}%
\bigstar & \bigstar & B\left[  \Sigma~0\right] \\
\bigstar & \bigstar & D\left[  \Sigma~0\right] \\
0 & 0 & 0_{m_{1}}%
\end{array}
\right]
\]
and%
\[
\lambda F=\left[
\begin{array}
[c]{ccc}%
\bigstar & \bigstar & 0\\
\bigstar & \bigstar & \lambda\left[  \Sigma~0\right] \\
0 & 0 & 0_{m_{1}}%
\end{array}
\right]  \text{.}%
\]
Since the block $\left[  \Sigma~0\right]  $ has full row rank, we conclude
that $B=0$ and $D=\lambda I_{m_{2}}$. Moreover, $A^{\prime}$ must be
nonsingular because $[A^{\prime}~B]$ has full row rank. Now examine%
\[
F^{2}=\left[
\begin{array}
[c]{ccc}%
A^{\prime} & 0 & 0\\
C & \lambda I_{m_{2}} & \left[  \Sigma~0\right] \\
0 & 0 & 0_{m_{1}}%
\end{array}
\right]  ^{2}=\left[
\begin{array}
[c]{ccc}%
\left(  A^{\prime}\right)  ^{2} & 0 & 0\\
CA^{\prime}+\lambda C & \lambda^{2}I_{m_{2}} & \lambda\left[  \Sigma~0\right]
\\
0 & 0 & 0_{m_{1}}%
\end{array}
\right]
\]
and%
\[
\lambda F=\left[
\begin{array}
[c]{ccc}%
\lambda A^{\prime} & 0 & 0\\
\lambda C & \lambda^{2}I_{m_{2}} & \lambda\left[  \Sigma~0\right] \\
0 & 0 & 0_{m_{1}}%
\end{array}
\right]  \text{,}%
\]
so that $\left(  A^{\prime}\right)  ^{2}=\lambda A^{\prime}$ (and $A^{\prime}$
is nonsingular), and $CA^{\prime}+\lambda C=\lambda C$. The first of these
identities tells us that $A^{\prime}=\lambda I_{n-m_{1}-m_{2}}$, and the
second tells us that $C=0$. Thus, $A$ is unitarily *congruent to%
\[
\lambda I_{n-m_{2}-m_{1}}\oplus\left[
\begin{array}
[c]{cc}%
\lambda I_{m_{2}} & \left[  \Sigma~0\right] \\
0 & 0_{m_{1}}%
\end{array}
\right]  \text{,\quad}\Sigma=\operatorname{diag}(\sigma_{1},\ldots
,\sigma_{m_{2}})\text{, all }\sigma_{i}>0\text{,}%
\]
which is unitarily *congruent (permutation similar) to%
\[
\lambda I_{n-m_{2}-m_{1}}\oplus%
{\displaystyle\bigoplus\limits_{i=1}^{m_{2}}}
\left[
\begin{array}
[c]{cc}%
\lambda & \sigma_{i}\\
0 & 0
\end{array}
\right]  \oplus0_{m_{1}-m_{2}}\text{.}%
\]
(b) The singular values of the 2-by-2 blocks are $0$ and $\tau_{i}%
=\sqrt{\left\vert \lambda\right\vert ^{2}+\sigma_{i}^{2}}>\left\vert
\lambda\right\vert $, so $\sigma_{i}=\sqrt{\tau_{i}^{2}-\left\vert
\lambda\right\vert ^{2}}$. \medskip%
\newline
(c) $A-\lambda I$ is unitarily *congruent to%
\[
0_{n-m_{2}-m_{1}}\oplus%
{\displaystyle\bigoplus\limits_{i=1}^{m_{2}}}
\left[
\begin{array}
[c]{cc}%
0 & \sigma_{i}\\
0 & -\lambda
\end{array}
\right]  \oplus(-\lambda)I_{m_{1}-m_{2}}\text{,}%
\]
so its singular values are: $\tau_{1},\ldots,\tau_{m_{2}}$, $|\lambda|$ with
multiplicity $m_{1}-m_{2}$, and $0$ with multiplicity $n-m_{1}$.\hfill
\end{proof}

\medskip Let $q(t)=(t-\lambda_{1})(t-\lambda_{2})$ be a given quadratic
polynomial (possibly $\lambda_{1}=\lambda_{2}$). If $q(t)$ is the minimal
polynomial of a given $A\in M_{n}$, then $q(A)=0$, $A-\lambda_{1}I$ is a
$\lambda$-projection with $\lambda:=\lambda_{2}-\lambda_{1}$, and $A$ is not a
scalar matrix. Theorem \ref{LambdaProjection} gives a canonical form to which
$A-\lambda_{1}I$ (and hence $A$ itself) can be reduced by unitary *congruence.

\begin{corollary}
\label{QuadraticMinimalPoly}Suppose the minimal polynomial of a given $A\in
M_{n}$ has degree two, and suppose that $\lambda_{1},\lambda_{2}$ are the
eigenvalues of $A$ with respective multiplicities $d$ and $n-d$; if
$\lambda_{1}=\lambda_{2}$, let $d=n$. Suppose that $|\lambda_{1}|\geq
|\lambda_{2}|$ and let $\sigma_{1},\ldots,\sigma_{m}$ be the singular values
of $A $ that are strictly greater than $|\lambda_{1}|$ ($m=0$ is possible).
Then: (a) $A$ is unitarily *congruent to%
\begin{equation}
\lambda_{1}I_{n-d-m}\oplus%
{\displaystyle\bigoplus\limits_{i=1}^{m}}
\left[
\begin{array}
[c]{cc}%
\lambda_{1} & \gamma_{i}\\
0 & \lambda_{2}%
\end{array}
\right]  \oplus\lambda_{2}I_{d-m}\text{,}\label{QuadraticCanonical}%
\end{equation}
in which each
\[
\gamma_{i}=\sqrt{\sigma_{i}^{2}+|\lambda_{1}\lambda_{2}|^{2}\sigma_{i}%
^{-2}-|\lambda_{1}|^{2}-|\lambda_{2}|^{2}}>0\text{.}%
\]
The direct sum (\ref{QuadraticCanonical}) is uniquely determined, up to
permutation of summands, by the eigenvalues and singular values of $A$. The
singular values of $A$ are $\sigma_{1},\ldots,\sigma_{m}$, $|\lambda
_{1}\lambda_{2}|\sigma_{1}^{-1},\ldots,|\lambda_{1}\lambda_{2}|\sigma_{m}%
^{-1}$, $|\lambda_{1}|$ with multiplicity $n-d-m$, and $|\lambda_{2}|$ with
multiplicity $d-m$. \medskip%
\newline
(b) Two square complex matrices that have quadratic minimal polynomials are
unitarily *congruent if and only if they have the same eigenvalues and the
same singular values.
\end{corollary}

\begin{proof}
If $A$ is singular, then $\lambda_{2}=0$, $A$ is a $\lambda_{1}$-projection,
and the validity of the assertions of the corollary is ensured by Theorem
\ref{LambdaProjection}.

Now assume that $A$ is nonsingular. The hypotheses ensure that $A-\lambda
_{1}I$ is singular and nonzero, and that it is a $\lambda$-projection with
$\lambda:=\lambda_{2}-\lambda_{1}$. It is therefore unitarily *congruent to a
direct sum of the form (\ref{ProjectorBlocks}) with $m_{1}=d$ and $m=m_{2}$;
after a translation by $\lambda_{1}I$, we find that $A$ is unitarily
*congruent to a direct sum of the form (\ref{QuadraticCanonical}) in which
each $\gamma_{i}=(\tau_{i}^{2}-\left\vert \lambda\right\vert ^{2})^{1/2}>0$.
Therefore, $A$ has $n-d-m$ singular values equal to $|\lambda_{1}|$ and $d-m$
singular values equal to $|\lambda_{2}|$. In addition, corresponding to each
2-by-2 block in (\ref{QuadraticCanonical}) is a pair of singular values
$\left(  \sigma_{i},\rho_{i}\right)  $ of $A$ such that
\begin{equation}
\sigma_{i}\geq\rho_{i}>0\text{ and }\sigma_{i}\rho_{i}=|\lambda_{1}\lambda
_{2}|\label{constraints}%
\end{equation}
for each $i=1,\ldots,m$. Since the spectral norm always dominates the spectral
radius, we have $\sigma_{i}\geq|\lambda_{1}|$ for each $i=1,\ldots,m$;
calculating the Frobenius norm tells us that
\[
\sigma_{i}^{2}+\rho_{i}^{2}=|\lambda_{1}|^{2}+|\lambda_{2}|^{2}+\gamma_{i}%
^{2}\text{.}%
\]
If $\sigma_{i}=|\lambda_{1}|$ then (\ref{constraints}) ensures that $\rho
_{i}=|\lambda_{2}|$, which is impossible since $\gamma_{i}>0$. Thus, $A$ has
$m$ singular values $\sigma_{1},\ldots,\sigma_{m}$ that are strictly greater
than $|\lambda_{1}|$, and $m$ corresponding singular values $\rho_{1}%
,\ldots,\rho_{m}$ that are strictly less than $|\lambda_{1}|$; each pair
$\left(  \sigma_{i},\rho_{i}\right)  $ satisfies (\ref{constraints}). Thus,
the parameters $\gamma_{i}$ in (\ref{QuadraticCanonical}) satisfy%
\[
\gamma_{i}^{2}=\sigma_{i}^{2}+\rho_{i}^{2}-|\lambda_{1}|^{2}-|\lambda_{2}%
|^{2}=\sigma_{i}^{2}+|\lambda_{1}\lambda_{2}|^{2}\sigma_{i}^{-2}-|\lambda
_{1}|^{2}-|\lambda_{2}|^{2}\text{.}%
\]

If two complex matrices of the same size have quadratic minimal polynomials,
and if they have the same eigenvalues and singular values, then each is
unitarily *congruent to a direct sum of the form (\ref{QuadraticCanonical}) in
which the parameters $\lambda_{1},\lambda_{2},d,m$, and $\{\gamma_{1}%
,\ldots,\gamma_{m}\}$ are the same; the two direct sums must be the same up to
permutation of their summands. Conversely, \textit{any} two unitarily
*congruent matrices have the same eigenvalues and singular values.\hfill
\end{proof}

Let $p(t)=t^{2}-2at+b^{2}$ be a given monic polynomial of degree two.
Corollary \ref{QuadraticMinimalPoly} tells us that if $p(A)=0$, then $A$ is
unitarily *congruent to a direct sum of certain special 1-by-1 and 2-by-2
blocks. We can draw a similar conclusion under the weaker hypothesis that
$p(A)$ is normal.

\begin{proof}
\begin{proposition}
\label{QuadraticPolynomialNormal}Let $A\in M_{n}$ and suppose there are
$a,b\in\mathbb{C}$ such that $N=A^{2}-2aA+bI$ is normal. Then $A$ is unitarily
*congruent to a direct sum of blocks, each of which is%
\[
\left[  \lambda\right]  \quad\text{or \quad}\left[
\begin{array}
[c]{cc}%
a+\nu & r\\
0 & a-\nu
\end{array}
\right]  ,\quad\lambda,\nu\in\mathbb{C},r\mathbf{\in}\mathbb{R},r>0,\text{ and
}\nu\in\mathcal{D}_{+}\text{.}%
\]

\end{proposition}
\end{proof}

\begin{proof}
A calculation reveals that $(A-aI)^{2}=N+(a^{2}-b)I$, which is normal. The
conclusion follows from applying Theorem \ref{squared normal canonical 2} to
the squared normal matrix $A-aI$.\hfill
\end{proof}

\subsection{Characterizations}

Corollary \ref{ConjugateNormalUnitaryCongruenceCorollary} tells us that a
conjugate normal matrix is unitarily congruent to a direct sum of a zero
matrix and positive scalar multiples of real orthogonal matrices; such a
matrix is real and normal. The following theorem gives additional
characterizations of conjugate normal matrices.

\begin{theorem}
\label{CconjugateNormal}Let $A\in M_{n}$ and let $A=PU=UQ$ be left and right
polar decompositions. Let $\sigma_{1}>\sigma_{2}>\cdots>\sigma_{d}\geq0$ be
the ordered distinct singular values of $A$ with respective multiplicities
$n_{1},\ldots,n_{d}$ (if $A=0$ let $d=1$ and $\sigma_{1}=0$). Let
$A=\mathcal{S}+\mathcal{C}$, in which $\mathcal{S}=\left(  A+A^{T}\right)  /2$
is symmetric and $\mathcal{C}=\left(  A-A^{T}\right)  /2$ is skew symmetric.
The following are equivalent: \medskip%
\newline
(a) $\mathcal{S\bar{C}}=\mathcal{C\bar{S}}$. \medskip%
\newline
(b) $A$ is conjugate normal. \medskip%
\newline
(c) $Q=\bar{P}$, that is, $A=PU=U\bar{P}$. \medskip%
\newline
(d) $PA=A\bar{P}$. \medskip%
\newline
(e) There are unitary matrices $W_{1},\ldots,W_{d}$ with respective sizes
$n_{1},\ldots,n_{d}$ such that $A$ is unitarily congruent to
\begin{equation}
\sigma_{1}W_{1}\oplus\cdots\oplus\sigma_{d}W_{d}.\label{e1e}%
\end{equation}
(f) There are real orthogonal matrices $Q_{1},\ldots,Q_{d}$ with respective
sizes $n_{1},\ldots,n_{d}$ such that $A$ is unitarily congruent to the real
normal matrix
\begin{equation}
\sigma_{1}Q_{1}\oplus\cdots\oplus\sigma_{d}Q_{d}.\label{e2e}%
\end{equation}

\end{theorem}

\begin{proof}
(a) $\Leftrightarrow$ (b): Compute%
\begin{align*}
A^{\ast}A  & =\left(  \mathcal{\bar{S}}-\mathcal{\bar{C}}\right)  \left(
\mathcal{S}+\mathcal{C}\right)  =\mathcal{\bar{S}S+\bar{S}C-\bar{C}S-\bar{C}%
C}\\
\overline{AA^{\ast}}  & =\left(  \mathcal{\bar{S}}+\mathcal{\bar{C}}\right)
\left(  \mathcal{S}-\mathcal{C}\right)  =\mathcal{\bar{S}S-\bar{S}C+\bar
{C}S-\bar{C}C}\\
A^{\ast}A-\overline{AA^{\ast}}  & =2\left(  \mathcal{\bar{S}C-\bar{C}%
S}\right)  \text{.}%
\end{align*}
Thus, $A^{\ast}A=\overline{AA^{\ast}}$ if and only if $\mathcal{\bar{S}%
C}=\mathcal{\bar{C}S}$.\medskip%
\newline
(b) $\Rightarrow$ (c): If $p(t)$ is any polynomial such that $p(\sigma_{i}%
^{2})=\sigma_{i}$ for each $i=1,...,d$, then $Q=p\left(  A^{\ast}A\right)  $
and $P=p(AA^{\ast})$. If $A^{\ast}A=\overline{AA^{\ast}}$ then
\[
Q=p\left(  A^{\ast}A\right)  =p(\overline{AA^{\ast}})=p\left(  \left(
AA^{\ast}\right)  ^{T}\right)  =p\left(  AA^{\ast}\right)  ^{T}=P^{T}=\bar
{P}\text{.}%
\]
(c) $\Rightarrow$ (d): $A\bar{P}=P(U\bar{P})=P(PU)=PA$.\medskip%
\newline
(d) $\Rightarrow$ (c): Let $P=V\Lambda V^{\ast}$ in which $V$ is unitary and
$\Lambda$ is nonnegative diagonal. Let $W=V^{\ast}U\bar{V}$. Then%
\[
A\bar{P}=PU\bar{P}=(V\Lambda V^{\ast})U(\bar{V}\Lambda V^{T})=V(\Lambda
W\Lambda)V^{T}%
\]
and%
\[
PA=P^{2}U=V\Lambda^{2}V^{\ast}U=V(\Lambda^{2}W)V^{T}\text{,}%
\]
so $\Lambda W\Lambda=\Lambda^{2}W$. Lemma \ref{Commutivity Implication}
ensures that $\Lambda W=W\Lambda$, so%
\[
PU=V\Lambda V^{\ast}U\bar{V}V^{T}=V\Lambda WV^{T}=VW\Lambda V^{T}=U\bar
{V}\Lambda V^{T}=U\bar{P}\text{.}%
\]
\medskip(c) $\Rightarrow$ (e): Suppose $P=V\Lambda V^{\ast}$, in which
$\Lambda=\sigma_{1}I_{n_{1}}\oplus\cdots\oplus\sigma_{d}I_{n_{d}}$ and $V$ is
unitary. If $Q=\bar{P}$ then%
\[
A=PU=V\Lambda V^{\ast}U=U\bar{V}\Lambda V^{T}=U\bar{P}=UQ=A
\]
and hence
\[
\Lambda\left(  V^{\ast}U\bar{V}\right)  =\left(  V^{\ast}U\bar{V}\right)
\Lambda\text{,}%
\]
which implies that the unitary matrix $V^{\ast}U\bar{V}=W_{1}\oplus
\cdots\oplus W_{d}$ is block diagonal; each $W_{i}$ is unitary and has size
$n_{i}$. Thus,
\[
U=V\left(  W_{1}\oplus\cdots\oplus W_{d}\right)  V^{T}%
\]
and%
\begin{align*}
A  & =PU=V\Lambda V^{\ast}U=V\Lambda V^{\ast}V\left(  W_{1}\oplus\cdots\oplus
W_{d}\right)  V^{T}\\
& =V\left(  \sigma_{1}W_{1}\oplus\cdots\oplus\sigma_{d}W_{d}\right)
V^{T}\text{.}%
\end{align*}
(e) $\Rightarrow$ (f): Corollary \ref{UnitarySpecialCase} ensures that each
$W_{j}$ in (\ref{e1e}) is unitarily congruent to a real orthogonal
matrix.\medskip%
\newline
(f) $\Rightarrow$ (a): Let $Z$ denote the direct sum (\ref{e2e}) and suppose
$A=UZU^{T}$ for some unitary $U$. Then $\mathcal{S}=%
\frac12
U(Z+Z^{T})U^{T}$ and $\mathcal{C}=%
\frac12
U(Z-Z^{T})U^{T}$, so it suffices to show that $Z$ commutes with $Z^{T}$. But
each $Q_{i}$ is real orthogonal, so%
\[
ZZ^{T}=\sigma_{1}^{2}Q_{1}Q_{1}^{T}\oplus\cdots\oplus\sigma_{d}^{2}Q_{d}%
Q_{d}^{T}=\sigma_{1}^{2}I_{n_{1}}\oplus\cdots\oplus\sigma_{d}^{2}I_{n_{d}%
}=Z^{T}Z\text{.}%
\]

\hfill
\end{proof}

\medskip

For normal matrices, an analog of Theorem \ref{CconjugateNormal} is the
following set of equivalent statements: \medskip%
\newline
(a) $HK=KH$, in which $H=(A+A^{\ast})/2$ and $K=(A-A^{\ast})/(2i)$. \medskip%
\newline
(b) $A$ is normal. \medskip%
\newline
(c) $Q=P$, that is, $A=PU=UP$. \medskip%
\newline
(d) $PA=AP$. \medskip%
\newline
(e) $A$ is unitarily *congruent to a direct sum of the form (\ref{e1e}), in
which $\sigma_{1}>\cdots>\sigma_{d}\geq0$ are the distinct singular values of
$A$ and $W_{1},\ldots,W_{d}$ are unitary.

\medskip The following theorem about conjugate normal matrices is an analog of
a known result about *congruence of normal matrices \cite{Ikramov} (and, more
generally, about unitoid matrices \cite[p. 289]{JF Sylvester}).

\begin{theorem}
(a) A nonsingular complex matrix is congruent to a conjugate normal matrix if
and only if it is congruent to a unitary matrix. \medskip%
\newline
(b) A singular complex matrix is congruent to a conjugate normal matrix if and
only if it is congruent to a direct sum of a unitary matrix and a zero
matrix.\medskip%
\newline
(c) Each conjugate normal matrix $A\in M_{n}$ is congruent to a direct sum,
uniquely determined up to permutation of summands, of the form%
\begin{equation}
I_{r-2q}\oplus%
{\displaystyle\bigoplus\limits_{j=1}^{q}}
\left[
\begin{array}
[c]{cc}%
0 & 1\\
e^{i\theta_{j}} & 0
\end{array}
\right]  \oplus0_{n-r}\text{,\quad}0<\theta_{j}\leq\pi\text{,}%
\label{conjugateNormalCanonical}%
\end{equation}
in which $r=\operatorname{rank}A$ and there is one block $H_{2}(e^{i\theta
_{j}})$ corresponding to each eigenvalue of $\bar{A}A$ that lies on the open
ray $\{te^{i\theta_{j}}:t>0\}$. The summand $I_{r-2q}$ corresponds to the
$r-2q$ positive eigenvalues of $\bar{A}A$. \medskip%
\newline
(d) Two conjugate normal matrices $A$ and $B$ of the same size are congruent
if and only if for each $\theta\in\lbrack0,\pi]$, $\bar{A}A$ and $\bar{B}B$
have the same number of eigenvalues on each open ray $\{te^{i\theta}:t>0\}$.
\end{theorem}

\begin{proof}
Only assertion (d) requires comment. If $A$ is conjugate normal and
nonsingular, the decomposition (\ref{e1e}) ensures that $\bar{A}A$ is
unitarily similar to (and hence has the same eigenvalues as)%
\[
\mathcal{W}=\sigma_{1}^{2}\overline{W_{1}}W_{1}\oplus\cdots\oplus\sigma
_{d}^{2}\overline{W_{d}}W_{d}\text{.}%
\]
Of course, $\mathcal{W}$ and the unitary matrix%
\[
W_{1}\overline{W_{1}}\oplus\cdots\oplus W_{d}\overline{W_{d}}%
\]
have the same number of eigenvalues on each open ray $\{te^{i\theta}:t>0\}$;
this number is the same as the number of blocks $H_{2}(e^{i\theta_{j}})$ in
the direct sum (\ref{conjugateNormalCanonical}). The argument is similar if
$A$ is singular; just omit the last direct summand $\sigma_{d}W_{d}$.\hfill
\end{proof}

\medskip

There is an analog of Theorem \ref{CconjugateNormal} for congruence normal matrices.

\begin{theorem}
\label{CongruenceNormalCharacterization}Let $A\in M_{n}$ and let $A=PU=UQ$ be
left and right polar decompositions. Let $A=\mathcal{S}+\mathcal{C}$, in which
$\mathcal{S}=(A+A^{T})/2$ is symmetric and $\mathcal{C}=(A-A^{T})/2$ is skew
symmetric. The following are equivalent: \medskip%
\newline
(a) $\mathcal{\bar{S}S}+\mathcal{\bar{C}C}$ commutes with $\mathcal{\bar{S}%
C}+\mathcal{\bar{C}S}$. \medskip%
\newline
(b) $A$ is congruence normal. \medskip%
\newline
(c) $A\bar{P}=\bar{Q}A$. \medskip%
\newline
(d) $\{\bar{P},Q,\bar{U}U\}$ is a commuting family.
\end{theorem}

\begin{proof}
(a) $\Leftrightarrow$ (b): A computation reveals that the Hermitian part of
$\bar{A}A$ is $\mathcal{\bar{S}S}+\mathcal{\bar{C}C}$, while the
skew-Hermitian part is $\mathcal{\bar{S}C}+\mathcal{\bar{C}S}$. Of course,
$\bar{A}A$ is normal if and only if its Hermitian and skew-Hermitian parts
commute. \medskip%
\newline
(b) $\Leftrightarrow$ (c): Theorem \ref{NonsingularEquivalence} tells us that
if $A$ is congruence normal then $A\left(  \bar{P}\right)  ^{2}=\left(
\bar{Q}\right)  ^{2}A$, which is the same as $A\left(  P^{2}\right)
^{T}=\left(  Q^{2}\right)  ^{T}A$, which implies that $Ap(P^{2})^{T}%
=p(Q^{2})^{T}A$ for any polynomial $p(t)$. Choose $p(t)$ such that
$p(t)=+\sqrt{t}$ on the spectrum of $P$ (and hence also on the spectrum of
$Q$), and conclude that $AP^{T}=Q^{T}A$, or $A\bar{P}=\bar{Q}A$. The converse
implication is immediate:
\[
A\bar{P}=\bar{Q}A\Rightarrow A\left(  \bar{P}\right)  ^{2}=\left(  \bar
{Q}\right)  ^{2}A\text{.}%
\]
(b) $\Rightarrow$ (d): Suppose $V$ is unitary, let $\mathcal{A}:=VAV^{T}$, and
consider the factors of the left and right polar decompositions $\mathcal{A}%
=\mathcal{PU}=\mathcal{UQ}$. One checks that $\mathcal{P}=VPV^{\ast}$,
$\mathcal{Q}=\bar{V}QV^{T}$, and $\mathcal{U}=VUV^{T}$. Moreover, $\{\bar
{P},Q,\bar{U}U\}$ is a commuting family if and only if $\{\mathcal{\bar{P}%
},\mathcal{Q},\overline{\mathcal{U}}\mathcal{U}\}$ is a commuting family.
Thus, if $A$ is congruence normal, there is no lack of generality to assume
that it is a direct sum of blocks of the form (\ref{gcon0}). Blocks of the
first type in (\ref{gcon0}) are 1-by-1, so commutation is trivial. For blocks
of the second type, the polar factors are%
\[
P=\tau\left[
\begin{array}
[c]{cc}%
1 & 0\\
0 & |\mu|
\end{array}
\right]  \text{, }Q=\tau\left[
\begin{array}
[c]{cc}%
|\mu| & 0\\
0 & 1
\end{array}
\right]  \text{, and }U=\left[
\begin{array}
[c]{cc}%
0 & 1\\
e^{i\theta} & 0
\end{array}
\right]  \text{, so }\bar{U}U=\left[
\begin{array}
[c]{cc}%
e^{i\theta} & 0\\
0 & e^{-i\theta}%
\end{array}
\right]  \text{.}%
\]
For both types of blocks, $\{\bar{P},Q,\bar{U}U\}$ is a diagonal family, so it
is a commuting family. \medskip%
\newline
(d) $\Rightarrow$ (b): If $\{\bar{P},Q,\bar{U}U\}$ is a commuting family, then%
\begin{align*}
\bar{A}A  & =\bar{P}\bar{U}UQ=\bar{P}\left(  \bar{U}U\right)  Q=\left(
\bar{U}U\right)  \left(  \bar{P}Q\right) \\
& =\bar{P}\left(  \bar{U}U\right)  Q=\left(  \bar{P}Q\right)  \left(  \bar
{U}U\right)  \text{.}%
\end{align*}
Since $\bar{P}$ and $Q$ are commuting positive semidefinite Hermitian
matrices, $\bar{P}Q$ is positive semidefinite Hermitian. But $\bar{U}U$ is
unitary without further assumptions, so we have a polar decomposition of
$\bar{A}A$ in which the factors commute. This ensures that $\bar{A}A$ is
normal.\hfill
\end{proof}

\medskip A calculation reveals that if $\mathcal{S\bar{C}}=\mathcal{C\bar{S}}
$ then $\mathcal{\bar{S}S+\bar{C}C}$ commutes with $\mathcal{\bar{S}C+\bar
{C}S}$, and that if $PU=U\bar{P}$ then $\{\bar{P},Q,\bar{U}U\}$ is a commuting
family. Thus, Theorems \ref{CconjugateNormal} and
\ref{CongruenceNormalCharacterization} permit us to conclude (again) that
every conjugate normal matrix is congruence normal.

For squared normal matrices, an analog of Theorem
\ref{CongruenceNormalCharacterization} is the following set of equivalent
statements: \medskip%
\newline
(a) $H^{2}-K^{2}$ commutes with $HK+KH$, in which $H=(A+A^{\ast})/2$ and
$K=(A-A^{\ast})/(2i)$. \medskip%
\newline
(b) $A^{2}$ is normal. \medskip%
\newline
(c) $AP=QA$. \medskip%
\newline
(d) $\{P,Q,U^{2}\}$ is a commuting family.

\medskip Our final characterization links the parallel expositions we have
given for squared normality and congruence normality.

\begin{theorem}
\label{BarBlocks}Let $A\in M_{n}$ and let%
\begin{equation}
\mathcal{A}=\left[
\begin{array}
[c]{cc}%
0 & A\\
\bar{A} & 0
\end{array}
\right] \label{BarBlockDef}%
\end{equation}
Then:
\newline
(a) $A^{2}$ is normal if and only if $\mathcal{A}$ is congruence normal.
\newline
(b) $A$ is congruence normal if and only if $\mathcal{A}^{2}$ is normal.
\newline
(c) $A$ is normal if and only if $\mathcal{A}$ is conjugate normal.
\newline
(d) $A$ is conjugate normal if and only if $\mathcal{A}$ is normal.
\newline
(e) $\mathcal{A}\overline{\mathcal{A}}\mathcal{A}^{T}=\mathcal{A}^{T}%
\overline{\mathcal{A}}\mathcal{A}$ if and only if $A^{\ast}A^{2}=A^{2}A^{\ast
}$.
\newline
(f) $\mathcal{A}^{\ast}\mathcal{A}^{2}=\mathcal{A}^{2}\mathcal{A}^{\ast}$ if
and only if $A\bar{A}A^{T}=A^{T}\bar{A}A$. \medskip%
\newline
Now suppose that $A$ is nonsingular. Then:
\newline
(g) $A^{-T}A$ is normal (respectively, Hermitian, unitary) if and only if
$\mathcal{A}^{-\ast}\mathcal{A}$ is normal (respectively, Hermitian,
unitary).
\newline
(h) $A^{-\ast}A$ is normal (respectively, Hermitian, unitary) if and only if
$\mathcal{A}^{-T}\mathcal{A}$ is normal (respectively, Hermitian, unitary).
\end{theorem}

\begin{proof}
Each assertion follows from a computation. For example, (a) follows from%
\[
\overline{\mathcal{A}}\mathcal{A}=\left[
\begin{array}
[c]{cc}%
\bar{A}^{2} & 0\\
0 & A^{2}%
\end{array}
\right]  \text{,}%
\]
(g) follows from%
\[
\mathcal{A}^{-\ast}\mathcal{A}=\left[
\begin{array}
[c]{cc}%
\overline{A^{-T}A} & 0\\
0 & A^{-T}A
\end{array}
\right]  \text{,}%
\]
and (h) follows from%
\[
\mathcal{A}^{-T}\mathcal{A}=\left[
\begin{array}
[c]{cc}%
\overline{A^{-\ast}A} & 0\\
0 & A^{-\ast}A
\end{array}
\right]  \text{.}%
\]

\hspace*{\fill}
\end{proof}

Using Theorem \ref{BarBlocks}, we can show that Theorems
\ref{NonsingularEquivalence} and \ref{Nonsingular*equivalence} are actually
equivalent: First apply Theorem \ref{Nonsingular*equivalence} to
$\mathcal{A},$ which tells us that $\mathcal{A}^{2}$ is normal if and only if
$\mathcal{A}^{\ast}\mathcal{A}^{2}=\mathcal{A}^{2}\mathcal{A}^{\ast}$ if and
only if $\mathcal{A}^{-\ast}\mathcal{A}$ is normal (if $A$ is nonsingular).
Theorem \ref{BarBlocks} (b), (f), and (g) now ensure that $A$ is congruence
normal if and only if $A\bar{A}A^{T}=A^{T}\bar{A}A$ if and only if $A^{-T}A$
is normal (if $A$ is nonsingular). Thus, Theorem \ref{Nonsingular*equivalence}
implies Theorem \ref{NonsingularEquivalence}. The reverse implication follows
from applying Theorem \ref{NonsingularEquivalence} to $\mathcal{A}$ and using
Theorem \ref{BarBlocks} (a), (e), and (h).

A similar argument shows that the equivalence of Theorem
\ref{HermitianCosquaresTheorem} (a) and (b) (respectively, Theorem
\ref{NonsingularConjugateNormal} (a) and (b)) implies and is implied by the
equivalence of Theorem \ref{Hermitian*cosquaresTheorem} (a) and (b)
(respectively, Theorem \ref{Unitary*cosquaresTheorem} (a) and (b)).

\subsection{The classification problem for cubed normals is unitarily wild}

We have seen that there are simple canonical forms for squared normal matrices
under unitary *congruence, and also for congruence normal matrices under
unitary congruence. However, the situation for cubed normal matrices under
unitary *congruence (and for matrices $A$ such that $A\bar{A}A$ is normal,
under unitary congruence) is completely different; the classification problems
in these cases are very difficult.

A problem involving complex matrices is said to be \textit{unitarily wild }if
it contains the problem of classifying arbitrary square complex matrices under
unitary *congruence. Since the latter problem contains the problem of
classifying an arbitrary system of linear mappings on unitary spaces
\cite[Section 2.3]{VVSquiver}, it is reasonable to regard any unitarily wild
problem as hopeless (by analogy with nonunitary matrix problems that contain
the problem of classifying pairs of matrices under similarity \cite{B+S}).

Two lemmas are useful in showing that the unitary congruence classification
problems for (a) cubed normal matrices under unitary *congruence, and (b) for
matrices $A$ such that $A\bar{A}A$ is normal, are both unitarily wild.

\begin{lemma}
\label{Basic} Let $\lambda_{1},\ldots,\lambda_{d}$ be given distinct complex
numbers and let $F,F^{^{\prime}}\in M_{n}$ be given conformally partitioned
block upper triangular matrices%
\[
F=\left[
\begin{array}
[c]{cccc}%
\lambda_{1}I_{n_{1}} & F_{12} & \cdots & F_{1d}\\
& \lambda_{2}I_{n_{2}} & \cdots & F_{2d}\\
&  & \ddots & \vdots\\
0 &  &  & \lambda_{d}I_{n_{d}}%
\end{array}
\right]  ,\qquad F^{^{\prime}}=\left[
\begin{array}
[c]{cccc}%
\lambda_{1}I_{n_{1}} & F_{12}^{^{\prime}} & \cdots & F_{1d}^{^{\prime}}\\
& \lambda_{2}I_{n_{2}} & \cdots & F_{2d}^{^{\prime}}\\
&  & \ddots & \vdots\\
0 &  &  & \lambda_{d}I_{n_{d}}%
\end{array}
\right]
\]
in which $n_{1}+n_{2}+\cdots+n_{d}=n$. If $S\in M_{n}$ and $SF=F^{^{\prime}}%
S$, then $S$ is block upper triangular conformal to $F$. If, in addition, $S$
is normal, then $S$ is block diagonal conformal to $F$.
\end{lemma}

\begin{proof}
Partition $S=[S_{ij}]_{i,j=1}^{d}$ conformally to $F$. Compare corresponding
$(i,j)$ blocks of $SF$ and $F^{^{\prime}}S$ in the order $(d,1),(d,2),\ldots
,(d,d-1)$ to conclude that each of $S_{d,1},S_{d,2},\ldots,S_{d,d-1}$ is a
zero block. Then continue by comparing the blocks in positions
$(d-1,1),(d-1,2),\ldots,(d-1,d-2)$, etc. If $S$ is normal and block upper
triangular, then Lemma \ref{Zero Blocks}(a) ensures that it is block
diagonal.\hspace*{\fill}
\end{proof}

\begin{lemma}
\label{SVDunique}Let $\sigma_{1}>\sigma_{2}>\cdots>\sigma_{d}\geq0$ and
$\sigma_{1}^{^{\prime}}>\sigma_{2}^{^{\prime}}>\cdots>\sigma_{d}^{^{\prime}%
}\geq0$ be given nonnegative real numbers, let $D,D^{^{\prime}}\in M_{n}$ be
given conformally partitioned block diagonal matrices%
\[
D=\left[
\begin{array}
[c]{cccc}%
\sigma_{1}I_{n_{1}} &  &  & \\
& \sigma_{2}I_{n_{2}} &  & \\
&  & \ddots & \\
&  &  & \sigma_{d}I_{n_{d}}%
\end{array}
\right]  ,\qquad D^{^{\prime}}=\left[
\begin{array}
[c]{cccc}%
\sigma_{1}^{^{\prime}}I_{n_{1}} &  &  & \\
& \sigma_{2}^{^{\prime}}I_{n_{2}} &  & \\
&  & \ddots & \\
&  &  & \sigma_{d}^{^{\prime}}I_{n_{d}}%
\end{array}
\right]
\]
in which $n_{1}+n_{2}+\cdots+n_{d}=n$. If $U,V\in M_{n}$ are unitary and
$DU=VD^{^{\prime}}$, then $\sigma_{i}=\sigma_{i}^{^{\prime}}$ for each
$i=1,\ldots,d$, and there are unitary matrices $W_{1}\in M_{n_{1}}%
,\ldots,W_{d-1}\in M_{n_{d-1}}$ and $\tilde{U},\tilde{V}\in M_{n_{d}}$ such
that $U=W_{1}\oplus\cdots\oplus W_{d-1}\oplus\tilde{U}$ and $V=W_{1}%
\oplus\cdots\oplus W_{d-1}\oplus\tilde{V}$; if $\sigma_{d}>0$ then $\tilde
{U}=\tilde{V}$.
\end{lemma}

\begin{proof}
Let $A=DU=VD^{^{\prime}}$. The eigenvalues of $AA^{\ast}=D^{2}$ and $A^{\ast
}A=(D^{^{\prime}})^{2}$ are the same, so $D=D^{^{\prime}}$. Moreover,
\[
AA^{\ast}=\left(  DU\right)  \left(  DU\right)  ^{\ast}=D^{2}=VD^{2}V^{\ast}%
\]
and
\[
A^{\ast}A=U^{\ast}D^{2}U=\left(  VD\right)  ^{\ast}\left(  VD\right)  =D^{2},
\]
so $D^{2}$ commutes with both $U$ and $V$ and hence each of $U$ and $V$ is
block diagonal conformal to $D$. The identity $DU=VD$ ensures that the
diagonal blocks of $U$ and $V$ corresponding to each $\sigma_{i}>0$ are
equal.\hfill
\end{proof}

\begin{theorem}
\label{WildCube}The problem of classifying square complex matrices $A$ up to
unitary *congruence is unitarily wild in both of the following two cases:
\medskip%
\newline
(a) $A^{3}=0$. \medskip%
\newline
(b) $A$ is nonsingular and $A^{3}$ is normal.
\end{theorem}

\begin{proof}
(a) Let $F,F^{^{\prime}}\in M_{k}$ be given. Define%
\begin{equation}
A=\left[
\begin{array}
[c]{ccc}%
0_{k} & I_{k} & F\\
0_{k} & 0_{k} & I_{k}\\
0_{k} & 0_{k} & 0_{k}%
\end{array}
\right]  \quad\text{and\quad}A^{^{\prime}}=\left[
\begin{array}
[c]{ccc}%
0_{k} & I_{k} & F^{^{\prime}}\\
0_{k} & 0_{k} & I_{k}\\
0_{k} & 0_{k} & 0_{k}%
\end{array}
\right]  \text{,}\label{AandAprime}%
\end{equation}
so that $A^{3}=(A^{^{\prime}})^{3}=0$ for any choices of $F$ and $F^{^{\prime
}}$. Suppose $A$ and $A^{^{\prime}}$ are unitarily *congruent, that is,
suppose there is a unitary $U=[U_{ij}]_{i,j=1}^{3}\in M_{3k}$, partitioned
conformally to $A$, such that $AU=UA^{^{\prime}}$. Then $A^{2}U=(UA^{^{\prime
}})^{2}$; the $1,3$ blocks of $A^{2}$ and $(A^{^{\prime}})^{2}$ are $I_{k}$
and all their other blocks are $0_{k}$. Comparison of the first block rows and
the third block columns of both sides of $A^{2}U=(UA^{^{\prime}})^{2}%
$\ reveals that $U_{11}=U_{33}$ and $U_{31}=U_{32}=U_{21}=0_{k}$. It follows
that $U$ is block diagonal since it is normal and block upper triangular.
Comparison of the $1,3$ blocks of both sides of $AU=UA^{^{\prime}}$ shows that
$FU_{22}=U_{11}F^{^{\prime}}$; comparison of the $1,2$ blocks shows that
$U_{11}=U_{22}$. Thus, $A$ and $A^{^{\prime}}$ are unitarily *congruent if and
only if $F$ and $F^{^{\prime}}$ are unitarily *congruent. \medskip%
\newline
(b) Let $F,F^{^{\prime}}\in M_{k}$ be given. Define the two nonsingular
matrices%
\[
A=\left[
\begin{array}
[c]{ccc}%
\lambda I_{2k} & 0 & A_{13}\\
0 & \mu I_{3k} & G\\
0 & 0 & I_{3k}%
\end{array}
\right]  \quad\text{and\quad}A^{^{\prime}}=\left[
\begin{array}
[c]{ccc}%
\lambda I_{2k} & 0 & A_{13}^{^{\prime}}\\
0 & \mu I_{3k} & G\\
0 & 0 & I_{3k}%
\end{array}
\right]  \text{,}%
\]
in which $\lambda=(-1+i\sqrt{3})/2$ and $\mu=\bar{\lambda}$ are the two
distinct roots of $t^{2}+t+1=0$,%
\[
G=\left[
\begin{array}
[c]{ccc}%
3I_{k} & 0 & 0\\
0 & 2I_{k} & 0\\
0 & 0 & I_{k}%
\end{array}
\right]  \text{,\quad}A_{13}=\left[
\begin{array}
[c]{ccc}%
I_{k} & I_{k} & F\\
0_{k} & I_{k} & I_{k}%
\end{array}
\right]  \text{,\quad and }A_{13}^{^{\prime}}=\left[
\begin{array}
[c]{ccc}%
I_{k} & I_{k} & F^{^{\prime}}\\
0_{k} & I_{k} & I_{k}%
\end{array}
\right]  \text{.}%
\]
A computation reveals that $A^{3}=(A^{^{\prime}})^{3}=\lambda^{3}I_{2k}%
\oplus\mu^{3}I_{3k}\oplus I_{3k}$ is diagonal (and hence normal) for any
choices of $F$ and $F^{^{\prime}}$. Suppose $A$ and $A^{^{\prime}}$ are
unitarily *congruent, that is, suppose there is a unitary $U=[U_{ij}%
]_{i,j=1}^{3}\in M_{8k}$, partitioned conformally to $A$, such that
$AU=UA^{^{\prime}}$. Lemma \ref{Basic} ensures that $U$ is block diagonal.
Since
\[
\left(  AU\right)  _{23}=GU_{33}=U_{22}G=\left(  UA\right)  _{23}\text{,}%
\]
Lemma \ref{SVDunique} ensures that $U_{33}=U_{22}$ and that $U_{33}%
=V_{1}\oplus V_{2}\oplus V_{3}$ is block diagonal conformal to $G$. Partition
$U_{11}=[W_{ij}]_{i,j=1}^{2}$, in which $W_{11},W_{22}\in M_{k}$. Equating the
$1,3$ blocks of both sides of the identity $AU=UA^{^{\prime}}$ gives the
identity%
\[
A_{13}U_{33}=\left[
\begin{array}
[c]{ccc}%
I_{k}V_{1} & V_{2} & FV_{3}\\
0_{k} & V_{2} & V_{3}%
\end{array}
\right]  =\left[
\begin{array}
[c]{cc}%
W_{11} & W_{12}\\
W_{21} & W_{22}%
\end{array}
\right]  \left[
\begin{array}
[c]{ccc}%
I_{k} & I_{k} & F^{^{\prime}}\\
0_{k} & I_{k} & I_{k}%
\end{array}
\right]  =U_{11}A_{13}^{^{\prime}}\text{.}%
\]
Comparison of the $2,1$ blocks of both sides of this identity tells us that
$W_{21}=0$, so $W$ is block upper triangular and hence $W_{12}=0$ as well.
Comparison of the $2,2$ and $2,3$ blocks tells us that $V_{2}=W_{22}=V_{3}$
and comparison of the $1,2$ blocks tells us that $V_{2}=W_{11}$. Finally,
comparison of the $1,3$ blocks and using $V_{3}=W_{11}$ reveals that
$FV_{3}=V_{3}F^{^{\prime}}$, so $A$ and $A^{^{\prime}}$ are unitarily
*congruent if and only if $F$ and $F^{^{\prime}}$ are unitarily
*congruent.\hfill
\end{proof}

\begin{theorem}
(a) The problem of classifying square complex matrices $A$ such that $A\bar
{A}A=0$ up to unitary congruence contains the problem of classifying arbitrary
square matrices up to unitary congruence. \medskip%
\newline
(b) The problem of classifying square complex matrices up to unitary
congruence is unitarily wild.
\end{theorem}

\begin{proof}
(a) Suppose the matrices $A$ and $A^{^{\prime}}$ in (\ref{AandAprime}) are
unitarily congruent, that is, $AU=\bar{U}A^{^{\prime}}$ for some unitary
$U=[U_{ij}]_{i,j=1}^{3}$ that is partitioned conformally to $A$. Then%
\begin{equation}
\left[
\begin{array}
[c]{ccc}%
U_{21}+FU_{31} & U_{22}+FU_{32} & U_{23}+FU_{33}\\
U_{31} & U_{32} & U_{33}\\
0 & 0 & 0
\end{array}
\right]  =\left[
\begin{array}
[c]{ccc}%
0 & \bar{U}_{11} & \bar{U}_{11}F^{^{\prime}}+\bar{U}_{12}\\
0 & \bar{U}_{21} & \bar{U}_{21}F^{^{\prime}}+\bar{U}_{22}\\
0 & \bar{U}_{31} & \bar{U}_{31}F^{^{\prime}}+\bar{U}_{32}%
\end{array}
\right]  \text{.}\label{3by3}%
\end{equation}
Comparing the $2,1$ blocks of both sides of (\ref{3by3}) tells us that
$U_{31}=0$, and then comparing the $1,1$ blocks as well as the $3,3$ blocks
tells us that $U_{21}=U_{32}=0$. Since $U$ is block upper triangular and
normal, it is block diagonal. Comparing the $1,2$ blocks and the $2,3$ blocks
of (\ref{3by3}) now tells us that $\bar{U}_{11}=U_{22}=\bar{U}_{33}$, so
$U_{11}=U_{33}$. Finally, comparing the $1,3$ blocks reveals that
$FU_{11}=\bar{U}_{11}F^{^{\prime}}$, that is, $A$ and $A^{^{\prime}}$ are
unitarily congruent if and only if $F$ and $F^{^{\prime}}$ are unitarily
congruent. \medskip%
\newline
(b) Let $F,F^{^{\prime}}\in M_{k}$ be given and suppose that%
\[
A=\left[
\begin{array}
[c]{cccc}%
0_{k} & I_{k} & 0 & F\\
0 & 0_{k} & I_{k} & 0\\
0 & 0 & 0_{k} & I_{k}\\
0 & 0 & 0 & 0_{k}%
\end{array}
\right]  \quad\text{and\quad}A^{^{\prime}}=\left[
\begin{array}
[c]{cccc}%
0_{k} & I_{k} & 0 & F^{^{\prime}}\\
0 & 0_{k} & I_{k} & 0\\
0 & 0 & 0_{k} & I_{k}\\
0 & 0 & 0 & 0_{k}%
\end{array}
\right]
\]
are unitarily congruent, that is, $AU=\bar{U}A^{^{\prime}}$ for some unitary
$U=[U_{ij}]_{i,j=1}^{4}$ that is partitioned conformally to $A$. An adaptation
of the argument in part (a) shows that $U$ is block diagonal, $U_{22}=\bar
{U}_{11}$, $U_{33}=U_{11}$, and $U_{44}=\bar{U}_{11}$. Hence, $FU_{11}%
=U_{11}F^{^{\prime}}$. We conclude that $A$ and $A^{^{\prime}}$ are unitarily
congruent if and only if $F$ and $F^{^{\prime}}$ are unitarily
*congruent.\hfill
\end{proof}

\subsection{A bounded iteration\label{BoundedIteration}}

Suppose $A\in M_{n}$ is nonsingular and let $x_{0}\in\mathbb{C}^{n}$ be given.
Define $x_{1},x_{2},\ldots$ by%
\begin{equation}
A^{T}x_{k+1}+Ax_{k}=0\text{, }k=0,1,2,\ldots.\label{TransposeIteration}%
\end{equation}
Under what conditions on $A$ is the sequence $x_{1},x_{2},\ldots$ bounded for
all choices of $x_{0}$?

We have%
\[
x_{k+1}=-A^{-T}Ax_{k}=\cdots=(-1)^{k}\left(  A^{-T}A\right)  ^{k+1}%
x_{0}\text{,}%
\]
so boundedness of the solution sequence for all choices of $x_{0}$ requires
that no eigenvalue of the cosquare $A^{-T}A$ has modulus greater than 1.
Moreover, every Jordan block of any eigenvalue of modulus 1 must be 1-by-1.
Inspection of (\ref{JCFcosquare}) reveals that the Jordan Canonical Form of
$A^{-T}A$ must have the form (\ref{JCFdiagonalizableCosquare}), in which each
$\mu_{j}\neq1$ and $\left\vert \mu_{j}\right\vert =1$. Theorem
\ref{CongruenceCanonicalForms}(a) ensures that $A$ is congruent to a direct
sum of blocks of the two types%
\begin{equation}
\lbrack1]\text{ and }\left[
\begin{array}
[c]{cc}%
0 & 1\\
\mu & 0
\end{array}
\right]  \text{,\quad}\left\vert \mu\right\vert =1\neq\mu\text{.}%
\label{BoundedBlocks}%
\end{equation}
Corollary \ref{UnitarySpecialCase} ensures that the 2-by-2 blocks in
(\ref{BoundedBlocks}) may be replaced by 2-by-2 real orthogonal blocks
(\ref{RealBlock}) or by 2-by-2 Hermitian unitary blocks (\ref{HermitianBlock}).

Conversely, if $A=SUS^{T}$ for some nonsingular $S$ and unitary $U$, then%
\[
0=A^{T}x_{k+1}+Ax_{k}=SU^{T}S^{T}x_{k+1}+SUS^{T}x_{k}\text{,\quad}k=1,2,\ldots
\]
if and only if%
\[
\xi_{k+1}=\left(  -1\right)  ^{k}\left(  \bar{U}U\right)  ^{k+1}\xi
_{0}\text{,\quad}\xi_{k}:=S^{T}x_{k}\text{, }k=1,2,\ldots\text{.}%
\]
The sequence $\xi_{0},\xi_{1},...$ is bounded since $\bar{U}U$ is unitary. In
summary, we have the following

\begin{theorem}
\label{BoundedSequenceTheorem}Let $A\in M_{n}$ be nonsingular. The following
are equivalent: \medskip%
\newline
(a) The sequence $x_{1},x_{2}\ldots$ defined by%
\[
A^{T}x_{k+1}+Ax_{k}=0\text{,\quad}k=0,1,2,\ldots
\]
is bounded for each given $x_{0}\in\mathbb{C}^{n}$. \medskip%
\newline
(b) $A$ is congruent to a unitary matrix. \medskip%
\newline
(c) $A$ is congruent to a real orthogonal matrix. \medskip%
\newline
(d) $A$ is congruent to a Hermitian unitary matrix. \medskip%
\newline
(e) $A$ is congruent to a nonsingular conjugate normal matrix.
\end{theorem}

Parallel reasoning using Theorems \ref{CosquareCharacterize}(b) and
\ref{CongruenceCanonicalForms}(b) leads to similar conclusions about the
conjugate transpose version of (\ref{TransposeIteration}).

\begin{theorem}
Let $A\in M_{n}$ be nonsingular. The following are equivalent: \medskip%
\newline
(a) The sequence $x_{1},x_{2}\ldots$ defined by%
\[
A^{\ast}x_{k+1}+Ax_{k}=0\text{,\quad}k=0,1,2,\ldots
\]
is bounded for each given $x_{0}\in\mathbb{C}^{n}$. \medskip%
\newline
(b) $A$ is *congruent to a unitary matrix. \medskip%
\newline
(c) $A$ is diagonalizable by *congruence. \medskip%
\newline
(d) $A$ is *congruent to a nonsingular normal matrix.
\end{theorem}

\section{Some comments about previous work}

Lemma \ref{Fuglede}(a) is often called the Fuglede-Putnam Theorem.

The assertion in Corollary \ref{CriterionCorollary}(b) that two unitary
matrices are *congruent if and only if they are unitarily *congruent was
proved in \cite{JF Sylvester} with an elegant use of uniqueness of the polar decomposition.

The unitary congruence canonical form (\ref{ConinvolutionCongruence}) for a
coninvolutory matrix was proved in \cite[Theorem 1.5]{HM1}.

Wigner \cite{Wigner} obtained a unitary congruence canonical form
(\ref{UnitarySpecialCaseBlocks}) for unitary matrices in which the 2-by-2
blocks are the Hermitian unitary blocks (\ref{HermitianBlock}).

In \cite{Autonne}, Autonne used a careful study of uniqueness of the unitary
factors in the singular value decomposition to prove many basic results, for
example: a nonsingular complex symmetric matrix is diagonalizable under
unitary congruence; a complex normal matrix is unitarily similar to a diagonal
matrix; a real normal matrix is real orthogonally similar to a real block
diagonal matrix with 1-by-1 and 2-by-2 blocks, in which the latter are scalar
multiples of real orthogonal matrices; similar unitary matrices are unitarily
similar. Lemma \ref{SVDunique} is a special case of Autonne's uniqueness
theorem; for an exposition see \cite[Theorem 3.1.1$^{^{\prime}}$]{HJ2}.

Hua proved the canonical form (\ref{skewsymmetric}) for a nonsingular skew
symmetric matrix under unitary congruence in \cite[Theorem 7]{Hua}; Theorem 5
in the same paper is the corresponding canonical form for a nonsingular
symmetric matrix.

The first studies of conjugate normal and congruence normal matrices seem to
be \cite{VujicicI} and \cite{HerbutII}.

The canonical form (\ref{t*0}) for a squared normal matrix (and hence the
canonical form (\ref{g*0})) can be deduced from Lemma 2.2 of \cite{VVSquiver}).

Each squared normal matrix can be reduced to the form (\ref{t*0}) by employing
the key ideas in Littlewood's algorithm \cite{Littlewood} for reducing
matrices to canonical form by unitary similarity. An exposition of this
alternative approach to Theorem \ref{squared normal canonical 2}, as well as a
canonical form for real squared normal matrices under real orthogonal
congruences, is in \cite{FHS}.

D.\v{Z}.
\raisebox{1.5pt}{-}\!\!Dokovi\'{c}
proved the canonical form (\ref{ProjectorBlocks}) for ordinary
projections ($\lambda=1$) in \cite{Dj}; for a different proof see \cite[p.
46]{VVSquiver}. George and Ikramov \cite{GI} used D.\v{Z}.
\raisebox{1.5pt}{-}\!\!Dokovi\'{c}'s canonical
form to derive a decomposition of the form (\ref{Involution*Congruence2}) for
an involution; in addition, they used Specht's Criterion to prove Corollary
\ref{QuadraticMinimalPoly}(b). For an ordinary projection $P$, and without
employing any canonical form for $P$, Lewkowicz \cite{Lewkowicz} identified
all of the singular values of $P$ and $I-P$.

The block matrix (\ref{BarBlockDef}) and the characterization of conjugate
normal matrices in Theorem \ref{BarBlocks}(d) was studied in \cite[Proposition
2]{F+I}. The characterization of conjugate normal matrices via the criterion
in Theorem \ref{CconjugateNormal}(a) is in \cite[Proposition 3]{F+I}.

Theorem \ref{WildCube}(a) was proved in \cite[p. 45]{VVSquiver}.

In \cite{Ikr1997}, Ikramov proved that any matrix with a quadratic minimal
polynomial is unitarily *congruent to a direct sum of the form
(\ref{QuadraticCanonical}). His characterization of the positive parameters
$\gamma_{i}$ is different from ours: If $\lambda_{1}\neq\lambda_{2}$, he found
that $\gamma=\left\vert \lambda_{1}-\lambda_{2}\right\vert \tan\alpha$, in
which $\alpha$ is the angle between any pair of left and right $\lambda_{1}%
$-eigenvectors of the block%
\[
\left[
\begin{array}
[c]{cc}%
\lambda_{1} & \gamma_{i}\\
0 & \lambda_{2}%
\end{array}
\right]  \text{.}%
\]
This pleasant characterization fails if $\lambda_{1}=\lambda_{2}$; our
characterization (using eigenvalues and singular values) is valid for all
$\lambda_{1},\lambda_{2}$.

The authors learned about the bounded iteration problem in Section
\ref{BoundedIteration} from Leiba Rodman and Peter Lancaster, who solved it
using canonical pairs.

\end{document}